\def\Bbb#1{{\fam\msbfam\relax#1}}
\font\fivemsb=msbm5
\font\sevenmsb=msbm7
\font\tenmsb=msbm10
\newfam\msbfam
\textfont\msbfam=\tenmsb
\scriptfont\msbfam\sevenmsb
\scriptscriptfont\msbfam\fivemsb

\def\spc{{\Bbb C}}
\def\spr{{\Bbb R}}
\def\vf{\varphi}
\def\wt{\widetilde}
\def\wh{\widehat}
\def\dm{\diamond}
\def\a{\alpha}
\def\b{\beta}

\def\F{\Phi}

\def\Q{\Psi}

\def\g{\gamma}
\def\l{\lambda}
\def\m{\mu}
\def\n{\nu}
\def\r{\rho}
\def\s{\sigma}

\def\q{\psi}
\def\t{\theta}
\def\T{\Theta}
\def\x{\xi}
\def\y{\eta}
\def\ca{{\cal A}}

\def\ce{{\cal E}}
\def\cg{{\cal G}}

\def\cl{{\cal L}}

\def\la{\langle}
\def\ra{\rangle}
\def\lla{\langle\!\langle}
\def\rra{\rangle\!\rangle}

\documentstyle[verbatim]{amsart}

\theoremstyle{plain}
\newtheorem{theorem}{Theorem}[section]

\newtheorem{lemma}[theorem]{Lemma}
\newtheorem{proposition}[theorem]{Proposition}
\newtheorem{fact}{Fact}[section]

\theoremstyle{definition}
\newtheorem{definition}[theorem]{Definition}
\newtheorem{example}{Example}[section]

\numberwithin{equation}{section}

\begin{document}

\title[characterization theorems, intrinsic topologies]%
{General characterization theorems and intrinsic
     topologies in white noise analysis}
\author{Nobuhiro Asai}
\author{Izumi Kubo}
\author{Hui-Hsiung Kuo}
\address{Nobuhiro Asai: International Institute
for Advanced Studies\\Kizu, Kyoto 619-0225 \\ JAPAN}
\address{Izumi Kubo: Department of Mathematics \\
Graduate School of Science \\ Hiroshima University \\
Higashi-Hiroshima 739-8526 \\ JAPAN}
\address{Hui-Hsiung Kuo: Department of Mathematics\\
  Louisiana State University \\ Baton Rouge \\
LA 70803 \\ USA}

\maketitle

\begin{abstract}
Let $u$ be a positive continuous function on $[0, \infty)$
satisfying the conditions: (i) $\lim_{r\to\infty} r^{-1/2}
\log u(r)=\infty$, (ii) $\inf_{r\geq 0} u(r)=1$, (iii) 
$\lim_{r\to \infty}\break r^{-1}\log u(r)<\infty$, (iv) the
function $\log u(x^{2}),\,x\geq 0$, is convex. A Gel'fand
triple $[\ce]_{u} \subset (L^{2}) \subset [\ce]_{u}^{*}$
is constructed by making use of the Legendre transform of $u$
discussed in \cite {akk3}.  We
prove a characterization theorem for generalized functions
in $[\ce]_{u}^{*}$  
and also for test functions in $[\ce]_{u}$
in terms of their $S$-transforms under 
the same assumptions on $u$.
Moreover, we give an intrinsic topology for the space
$[\ce]_{u}$ of test functions and prove a characterization
theorem for measures. 
We briefly mention the relationship 
between our method and a recent work by Gannoun et al.\cite{ghor}.
Finally, conditions for carrying out white noise 
operator theory and Wick products are given.

\end{abstract}

\medskip
\section{Introduction} \label{sec:1} 

Let $\ce$ be a real topological vector space with topology 
generated by a sequence of inner product norms 
$\{|\cdot|_{p}\}_{p=0}^{\infty}$. We assume that $\ce$ is 
a complete metric space with respect to the metric
\begin{equation}
 d(\x, \y) = \sum_{p=0}^{\infty} {1\over 2^{p}}
 {|\x-\y|_{p} \over 1+|\x-\y|_{p}}, \qquad \x, \y\in
 \ce.  \notag
\end{equation}
In addition we assume the following conditions:
\begin{itemize}
\item[(a)] There exists a constant $0<\r<1$ such that
$|\cdot|_{0}\leq \r|\cdot|_{1}\leq \cdots\leq 
 \r^{p}|\cdot|_{p} \leq \cdots$.
\item[(b)] For any $p\geq 0$, there exists $q\geq p$ such 
that the inclusion $i_{q, p}: \ce_{q} \hookrightarrow 
\ce_{p}$ is a Hilbert-Schmidt operator. (Here $\ce_{p}$ 
is the completion of $\ce$ with respect to the norm 
$|\cdot|_{p}$.)
\end{itemize}

Let $\ce'$ and $\ce_{p}'$ denote the dual spaces of $\ce$
and $\ce_{p}$, respectively. We can use the Riesz
representation theorem to identify $\ce_{0}$ with its
dual space $\ce_{0}'$. Then we get the following continuous
inclusions:
\begin{equation}  
 \ce \hookrightarrow \ce_{p}  \hookrightarrow \ce_{0}
  \hookrightarrow \ce_{p}'  \hookrightarrow \ce',
  \qquad p\geq 0.   \notag
\end{equation}

The above condition (b) says that $\ce$ is a nuclear
space and so $\ce \subset \ce_{0} \subset \ce'$ is a 
Gel'fand triple. Let $\m$ be the Gaussian measure on $\ce'$
with the characteristic function given by
\begin{equation}
 \int_{\ce'} e^{i\la x, \x\ra}\,d\m(x) =
  e^{-{1\over 2}|\x|_{0}^{2}}, \qquad \x\in\ce. \notag
\end{equation}

The probability space $(\ce', \m)$ is often referred to as
a {\em white noise space}. For simplicity, we will use 
$(L^{2})$ to denote the complex Hilbert space $L^{2}(\m)$. 
By the Wiener-It\^o theorem, each $\vf\in (L^{2})$ can be 
uniquely represented by
\begin{equation} \label{eq:1-1}
\vf (x) = \sum_{n=0}^{\infty} I_{n}(f_{n})(x) = 
 \sum_{n=0}^{\infty} \la :\!x^{\otimes n}\!:, f_{n}\ra, 
 \qquad f_{n} \in \ce_{0}^{\wh\otimes n},
\end{equation}
where $I_{n}$ is the multiple Wiener integral of order $n$
and $:\!x^{\otimes n}\!:$ is the Wick tensor of $x\in\ce'$
(see page 33 in \cite{kuo96}.) Moreover, the $(L^{2})$-norm 
of $\vf$ is given by
\begin{equation} \label{eq:1-2}
\|\vf\|_{0} = \left(\sum_{n=0}^{\infty} n!|f_{n}|_{0}^{2}
 \right)^{1/2}.  
\end{equation}

Recently Cochran et al.~\cite{cks} have introduced a 
Gel'fand triple associated with the above Gel'fand triple
$\ce \subset \ce_{0} \subset \ce'$ and a sequence 
$\{\a(n)\}_{n=0}^{\infty}$ of positive real numbers 
satisfying the conditions:
\begin{itemize}
\item[(A1)] $\a(0)=1$ and $\inf_{n\geq 0} \a(n) 
 \s^{n} >0$ for some $\s\geq 1$.
\smallskip
\item[(A2)] $\lim_{n\to\infty} \left({\a(n) \over n!}
 \right)^{1/n} =0$.
\end{itemize}
Actually, a stronger condition $\inf_{n\geq 0} \a(n) >0$ 
in (A1) is assumed in \cite{cks}. However the above weaker 
condition for some $\s\geq 1$ in (A1) is strong enough to 
assure that the nuclear space $[\ce]_{\a}$ is a subspace 
of $(L^{2})$, a fact to be shown below. This weaker 
condition was first introduced in \cite{akk4}.

For $\vf\in (L^{2})$ being represented by Equation 
(\ref{eq:1-1}) and $p\geq 0$, define
\begin{equation} \label{eq:1-3}
 \|\vf\|_{p, \a} = \left(\sum_{n=0}^{\infty} n! \a(n)
   |f_{n}|_{p}^{2}\right)^{1/2}.
\end{equation}
Let $[\ce_{p}]_{\a} = \{\vf\in (L^{2}); \,\|\vf\|_{p, \a}
< \infty\}$. Define the space $[\ce]_{\a}$ of {\em test
functions} on $\ce'$ to be the projective limit of 
$\{[\ce_{p}]_{\a}; \,p\geq 0\}$. The dual space 
$[\ce]_{\a}^{*}$ of $[\ce]_{\a}$ is called the space of
{\em generalized functions} on $\ce'$.

By using conditions (a) and (A1) we see that
$$
 \sum_{n=0}^{\infty} n!\a(n) |f_{n}|_{p}^{2}
\geq
 \Big(\inf_{n\geq 0} \a(n) \s^{n}\Big)
 \sum_{n=0}^{\infty} n! |f_{n}|_{0}^{2}  
$$
for $p$ large enough such that 
$\s^{-1} \r^{-2p}\geq 1$.
This inequality, in view of Equations (\ref{eq:1-2}) and 
(\ref{eq:1-3}), implies that $[\ce_{p}]_{\a} \subset 
(L^{2})$ for all $p\geq (-2\log \r)^{-1} \log \s$. Hence
$[\ce]_{\a} \subset (L^{2})$ holds. By identifying 
$(L^{2})$ with its dual space we get the following 
continuous inclusions:
\begin{equation}  
 [\ce]_{\a} \hookrightarrow [\ce_{p}]_{\a}  
 \hookrightarrow (L^{2}) \hookrightarrow 
 [\ce_{p}]_{\a}^{*} \hookrightarrow [\ce]_{\a}^{*},  
 \qquad p \geq (-2\log \r)^{-1}\log \s,  \notag
\end{equation}
where $[\ce_{p}]_{\a}^{*}$ is the dual space of 
$[\ce_{p}]_{\a}$. Moreover, $[\ce]_{\a}$ is a nuclear 
space and so $[\ce]_{\a} \subset (L^{2}) \subset 
[\ce]_{\a}^{*}$ is a Gel'fand triple. 
This triple is refered to as a {\it CKS-space}.
Note that 
$[\ce]_{\a}^{*}=\cup_{p\geq 0} [\ce_{p}]_{\a}^{*}$ and 
for $p\geq (-2\log \r)^{-1}\log \s$, 
$\,[\ce_{p}]_{\a}^{*}$ is the completion of $(L^{2})$ 
with respect to the norm
\begin{equation} \label{eq:1-4}
 \|\vf\|_{-p, 1/\a} = \left(\sum_{n=0}^{\infty} {n! \over 
 \a(n)}\,|f_{n}|_{-p}^{2}\right)^{1/2}.
\end{equation}

For $\x$ belonging to the complexification $\ce_{c}$ of 
$\ce$, the renormalized exponential function $:\!e^{\la 
\cdot, \x\ra}\!:$ is defined by
\begin{equation}
 :\!e^{\la \cdot, \x\ra}\!: \, = \sum_{n=0}^{\infty}
 {1 \over n!}\,\la :\!\cdot^{\otimes n}\!:, 
 \x^{\otimes n} \ra  
 = e^{\la \cdot, \xi \ra - {1 \over 2}\la \xi, \xi \ra}
,\notag
\end{equation}
whose norm is evaluated as 
\begin{equation} \label{eq:1-5}
 \|\!:\!e^{\la \cdot, \x\ra}\!:\!\|_{p, \a}
 = G_{\a}\big(|\x|_{p}^{2}\big)^{1/2}, \quad p \geq 0, 
\end{equation}
where $G_{\a}$ is the exponential generating function
of the sequence $\{\a(n)\}$, i.e.,
\begin{equation}  \label{eq:1-6}
 G_{\a}(r) = \sum_{n=0}^{\infty} {\a(n) \over n!}\,
  r^{n}. 
\end{equation}
By condition (A2) $G_{\a}$ is an entire function. Hence
Equation (\ref{eq:1-5}) implies that $:\!e^{\la \cdot, 
\x\ra}\!:\,\in [\ce]_{\a}$ for all $\x\in\ce_{c}$. 

On the other hand, by Equation (\ref{eq:1-4}), we have
\begin{equation}\label{eq:1-7}
 \|\!:\!e^{\la \cdot, \x\ra}\!:\!\|_{-p,1/\a}
 = G_{1/\a}\big(|\x|_{-p}^{2}\big)^{1/2},  
\end{equation}
where $G_{1/\a}$ is the exponential generating function
of the sequence $\{1/\a(n)\}$, i.e.,
\begin{equation}  \label{eq:1-8}
 G_{1/\a}(r) = \sum_{n=0}^{\infty} {1 \over n!\a(n)}\,
  r^{n}. 
\end{equation}
It follows from condition (A1) that $G_{1/\a}$ is an
entire function.

For $\F\in [\ce]_{\a}^{*}$, its {\em $S$-transform}
$S\F$ is defined to be the function
\begin{equation}
 (S\F)(\x) = \lla \F, :\!e^{\la \cdot, \x\ra}\!:\rra,
 \qquad \x \in \ce_{c},  \notag
\end{equation}
where $\lla \cdot, \cdot \rra$ is the bilinear pairing
of $[\ce]_{\a}^{*}$ and $[\ce]_{\a}$. 

An important problem in white noise analysis is to
characterize generalized and test functions in terms of
their $S$-transforms. For this purpose we need the
following conditions:
\begin{itemize}
\item[(B1)] $\limsup_{n\to\infty} \left({n! \over \a(n)}
 \inf_{r>0} {G_{\a}(r) \over r^{n}}\right)^{1/n} < \infty$.
\smallskip
\item[($\wt{\text{B}}1$)] $\limsup_{n\to\infty} \left(
 n!\a(n) \inf_{r>0} {G_{1/\a}(r) \over r^{n}}\right)^{1/n}
  < \infty$.
\smallskip
\item[(B2)] The sequence $\g(n)={\a(n) \over n!}, n\geq 0,$ 
is log-concave, i.e., 
\begin{equation}
 \g(n) \g(n+2) \leq \g(n+1)^{2}, \qquad 
      \forall n\geq 0. \notag
\end{equation}
\item[($\wt{\text{B}}2$)] The sequence $\left\{{1 \over 
n!\a(n)}\right\}$ is log-concave.
\end{itemize}

It follows from Theorem 4.3 in \cite{cks} that condition
(B2) implies condition (B1). Similarly, condition
($\wt{\text{B}}2$) implies condition 
($\wt{\text{B}}1$), \cite{akk2}.
Characterization thorems are proved 
in \cite{cks} for generalized functions 
under (B2) and in \cite{akk2} for test functions 
under ($\wt{\text{B}}2$).
In fact, those conditions
can be replaced by weaker conditions.
We say that two sequences $\{a(n)\}$ and 
$\{b(n)\}$ of positive real numbers are {\em equivalent} 
if there exist $K_{1}, K_{2}, c_{1}, c_{2}>0$ such that 
for all $n$,
\begin{equation}
 K_{1} c_{1}^{n}a(n) \leq b(n) \leq  K_{2} c_{2}^{n}a(n).
    \notag
\end{equation}

Now, we state the weaker conditions for the sequence
$\{\a(n)\}$:

\medskip\noindent
$\underline{\text{Near-(B2)}}\;$ The sequence $\{\a(n)\}$ 
is equivalent to a sequence $\{\l(n)\}$ of positive real 
numbers such that $\{\l(n)/n!\}$ is log-concave.

\medskip\noindent
$\underline{\text{Near-($\wt{\text{B}}2$)}}\;$ The sequence 
$\{\a(n)\}$ is equivalent to a sequence $\{\l(n)\}$ of 
positive real numbers such that $\{{1\over n!\l(n)}\}$ is 
log-concave.

\medskip
As shown in Lemma \ref{lem:2-1} later,  near-(B2) 
and near-($\wt{\text{B}}2$) are equivalent to  
the conditions (B1) and ($\wt{\text{B}}1$), respectively. 
Then, we have the following theorems.

\begin{theorem} \label{thm:1-1}
If $F=S\F$ for $\F\in [\ce]_{\a}^{*}$, 
then $F$ satisfies the conditions:
\begin{itemize}
\item[(1)] For any $\x, \y \in \ce_{c}$, the function
$F(z\x+\y)$ is an entire function of $z\in\spc$.
\item[(2)] There exist constants $K, a, p\geq 0$ such
that
\begin{equation}
 |F(\x)| \leq K G_{\a}\big(a|\x|_{p}^{2}\big)^{1/2},
  \qquad \x \in\ce_{c}.  \notag
\end{equation}
\end{itemize}
Conversely, assume that condition near-(B2) holds and let 
$F\!: \ce_{c}\to \spc$ be a function satisfying 
conditions (1) and (2). Then $F=S\F$ for a unique 
generalized function $\F\in [\ce]_{\a}^{*}$.
\end{theorem}

\begin{theorem} \label{thm:1-2}
If $F=S\vf$ for $\vf\in [\ce]_{\a}$,
then $F$ satisfies the conditions:
\begin{itemize}
\item[(1)] For any $\x, \y \in \ce_{c}$, the function
$F(z\x+\y)$ is an entire function of $z\in\spc$.
\item[(2)] For any $a, p\geq 0$, there exists a constant
$K\geq 0$ such
that
\begin{equation} 
 |F(\x)| \leq K G_{1/\a}\big(a|\x|_{-p}^{2}\big)^{1/2},
   \qquad \x \in\ce_{c}.  \notag
\end{equation}
\end{itemize}
Conversely, assume that condition near-($\wt{\text{B}}2$) holds 
and let $F\!: \ce_{c}\to \spc$ be a function satisfying 
conditions (1) and (2). Then $F=S\vf$ for a unique test
function $\vf \in [\ce]_{\a}$.
\end{theorem}

Now, for a general sequence $\{\a(n)\}$, we cannot expect
to find the sums $G_{\a}$ in Equation (\ref{eq:1-6}) and
$G_{1/\a}$ in Equation (\ref{eq:1-8}) as elementary 
functions. Therefore, it is desirable to find elementary
functions to replace $G_{\a}$ and $G_{1/\a}$ in Theorems
 \ref{thm:1-1} and \ref{thm:1-2}. This leads to the 
concept of equivalence in the next definition.

\begin{definition}
Two positive functions $f$ and $g$ on $[0, \infty)$ are
called {\em equivalent} if there exist constants 
$c_{1}, c_{2}, a_{1}, a_{2} >0$ such that
\begin{equation}
 c_{1}f(a_{1}r) \leq g(r) \leq c_{2}f(a_{2}r), \qquad
   \forall\> r\in [0, \infty).  \notag
\end{equation}
\end{definition}

In order to find elementary functions that are equivalent
to $G_{\a}$ and $G_{1/\a}$ in Theorems \ref{thm:1-1} and
\ref{thm:1-2}, we have developed in \cite{akk3} the crucial 
mathematical machinery.

\begin{example} 
When $\a(n)=1$ for all $n$, condition (B2) is obviously 
satisfied and $G_{\a}(r)=G_{1/\a}(r)=e^{r}$. In this case 
defined by Hida-Kubo-Takenaka, 
Theorem \ref{thm:1-1} is due to Potthoff and Streit
\cite{ps}, while Theorem \ref{thm:1-2} is due to 
Kuo et al.~\cite{kps}.
\end{example}

\begin{example} 
When $\a(n)=(n!)^{\b},\, 0\leq \b<1$, condition (B2) is
 easily seen to be satisfied. The functions $G_{\a}$ and
 $G_{1/\a}$ are given respectively by
 \begin{equation} \label{eq:1-9}
 G^{(\b)}(r) = \sum_{n=0}^{\infty} {1\over (n!)^{1-\b}}\,
 r^{n}, \quad  G^{(-\b)}(r) = \sum_{n=0}^{\infty} 
  {1\over (n!)^{1+\b}}\, r^{n}.
\end{equation}
However we see that $G^{(\b)}$ and $G^{(-\b)}$ 
are equivalent to the functions 
\begin{equation} \label{eq:1-10}
 \!\!\!\!\! g_{\b}(r) = \exp\Big[(1-\b) r^{1\over 
 1-\b}\Big], \quad  g_{-\b}(r) = \exp\Big[(1+\b)
   r^{1\over 1+\b}\Big], 
\end{equation}
respectively.  Theorems \ref{thm:1-1} and \ref{thm:1-2} 
with the growth functions $g_{\b}$ and $g_{-\b}$, respectively,
are due to Kondratiev and Streit \cite{ks92} \cite{ks93} 
 (see also \cite{kuo96}.)
\end{example} 

\begin{example} 
When $\a(n)=b_{k}(n)$ (the Bell numbers of order $k$,)
condition (B1) is shown to be satisfied in \cite{cks}. However 
actually condition (B2) is satisfied \cite{akk1}. In this 
case, $G_{\a}(r) = \exp_{k}(x)/\exp_{k}(0)$ ($\exp_{k}(x)$
is the $k$-th iterated exponential function) and Theorem 
\ref{thm:1-1} is due to Cochran et al.~\cite{cks}. However, 
$G_{1/\a}$ is not an elementary function. However it is equivalent 
to the function
\begin{equation} \label{eq:1-11}
 w_{k}(r) = \exp\left[2\,\sqrt{r\log_{k-1}\sqrt{r}}\right],
\end{equation}
where the function $\log_{j}$ is defined inductively by
\begin{equation}
 \log_{1}(r)=\log(\max\{r, e\}), \quad \log_{j}(r) =
  \log_{1}(\log_{j-1}(r)), \quad j\geq 2. \notag
\end{equation}
\end{example}

\medskip
The first purpose of the present paper 
is to construct a CKS-space $[\ce]_{u} \subset (L^{2}) 
\subset [\ce]_{u}^{*}$ with a given growth function $u$ 
and to obtain the general characterization theorems by 
applying the results in \cite{akk3}. 
The second purpose is to give the intrinsic topology 
for $[\ce]_{u}$ and to show properties of the space 
relating to the features of $u$.
The basic idea is to start with a growth function $u$ 
and then apply the Legendre transform to get a Gel'fand 
triple $[\ce]_{u} \subset (L^{2}) \subset [\ce]_{u}^{*}$. 

We remark that Gannoun et al.~\cite{ghor} studied 
a similar Gel'fand triple consisting of spaces of 
entire functions governed by a convex function 
$\theta$ and its dual.  Their inclusions and duality 
are rather abstract.  
We will give comments about relationships between 
$u$ and $\theta$ in section \ref{sec:4}.  

Further we will discuss the characterization of 
measures in $[\ce]_{u}^{*}$ by an integrability 
condition in section 4.  
Ouerdiane and Rezgui \cite{or} showed 
the Bochner-Minlos theorem, which tells an integrability 
condition of measures in terms of growth order of 
characteristic functions. 

\medskip
\section{Legendre and dual Legendre transforms}
   \label{sec:2}

First we mention the following concepts which will be
frequently used. A positive function $f$ on $[0, \infty)$ 
is called
\begin{itemize}
\item[(1)] {\em log-concave} if the function $\log f$ is 
 concave on $[0, \infty)$;
\item[(2)] {\em log-convex} if the function $\log f$ is 
 convex on $[0, \infty)$;
\item[(3)] {\em (log, exp)-convex} if the function 
 $\log f(e^{x})$ is convex on $\spr$;
\item[(4)] {\em (log, $x^{k}$)-convex} if the function 
 $\log f(x^{k})$ is convex on $[0, \infty)$. Here $k>0$.
\end{itemize}   

\smallskip

It is easy to check that if $f$ is log-concave, then
the sequence $\{f(n)\}_{n=0}^{\infty}$ is log-concave. If
$f$ is increasing and (log, $x^{k}$)-convex for some $k>0$,
then $f$ is (log, exp)-convex (see Proposition 2.3(3) in
[3].)   Further, 
if $\{\b(n)/n!\}_{n=0}^{\infty}$ is log-concave and
$\b(0)=1$, then for any $n, m \geq 0$,
\begin{equation} \label{eq:add}
	\b(n+m) \leq {n+m \choose n} \b(n) \b(m)
	\leq 2^{n+m}\b(n)\b(m). 
\end{equation}

\smallskip
Let $C_{+, \log}$ denote the collection of all positive 
continuous functions $u$ on $[0, \infty)$ satisfying the 
condition:
\begin{equation}
 \lim_{r\to\infty} {\log u(r) \over \log r}=\infty. \notag
\end{equation}
The {\em Legendre transform} $\,\ell_{u}$ of $u \in 
C_{+, \log}$ is defined to be the function
\begin{equation}
 \ell_{u}(t) = \inf_{r>0} {u(r) \over r^{t}}, \qquad
  t\in [0, \infty). \notag
\end{equation}

\smallskip
 Let $C_{+, \,1/2}$ denote the collection of all positive 
 continuous functions $u$ on $[0, \infty)$ satisfying the 
 condition:
 \begin{equation}
  \lim_{r\to\infty} {\log u(r) \over \sqrt{r}}=\infty. \notag
 \end{equation}
The {\em dual Legendre transform} $\,u^{*}$ of $u \in 
C_{+, 1/2}$ is defined to be the function
\begin{equation}
 u^{*}(r) = \sup_{s\geq 0} {e^{2\sqrt{rs}} \over u(s)}, 
   \qquad  r\in [0, \infty). \notag
\end{equation}
Note that $C_{+, 1/2} \subset C_{+, \log}$. 
Assume that $u\in C_{+, \log}$ and  
$\lim_{n\to\infty} \ell_{u}(n)^{1/n} =0$. We define the 
$L$-{\em function} $\cl_{u}$ of $u$ by
\begin{equation} \label{eq:2-1}
  \cl_{u}(r) = \sum_{n=0}^{\infty} \ell_{u}(n) r^{n}.
\end{equation}

Now, let $u\in C_{+, 1/2}$ and assume that 
$\lim_{n\to\infty} \big(\ell_{u}(n)(n!)^{2}\big)^{-1/n}=0$.
We define the {\em $L^{\#}$-function} of $u$ by
\begin{equation} \label{eq:2-2}
  \cl_{u}^{\#}(r) = \sum_{n=0}^{\infty} {1 \over 
  \ell_{u}(n)(n!)^{2}}\, r^{n}. 
\end{equation}

For discussions in the rest of the paper, we will
need the following facts from papers \cite{akk3} \cite{akk4}.

\begin{fact} \label{fact:2-0}
The inclusion $C_{+,{1 \over 2}} \subset C_{+, \log}$ holds. 
If $u$ is increasing and $(\log,x^{2})$-convex, 
then $u$ is $(log,exp)$-convex. 
\end{fact}
 
\begin{fact} \label{fact:2-1}
Let $u\in C_{+, \log}$. Then the Legendre transform 
$\ell_{u}$ is log-concave. (Hence $\ell_{u}$ is 
continuous on $[0, \infty)$ and the sequence 
$\{\ell_{u}(n)\}_{n=0}^{\infty}$ is log-concave.)
\end{fact}
 
%

\begin{fact} \label{fact:2-3}
Let $u\in C_{+, \log}$ be (log, exp)-convex. Then 
\begin{itemize}
\smallskip
\item[(1)] $\ell_{u}(t)$ is decreasing for large $t$,
\smallskip
\item[(2)] $\lim_{t\to\infty} \ell_{u}(t)^{1/t} =0$, 
\smallskip
\item[(3)] $u(r) = \sup_{t\geq 0} \ell_{u}(t) r^{t}$ 
for all $r\geq 0$.
\end{itemize}
\end{fact}

\begin{fact} \label{fact:2-4}
Let $u\in C_{+, \log}$. We have the assertions:
\begin{itemize}
\smallskip
\item[(1)] $u$ is (log, $x^{k}$)-convex if and only if 
$\ell_{u}(t) t^{kt}$ is log-convex.
\smallskip
\item[(2)] If $u$ is (log, $x^{k}$)-convex, then for any 
integers $n, m\geq 0$,
\begin{equation}
 \ell_{u}(n) \ell_{u}(m) \leq \ell_{u}(0) 2^{k(n+m)}
  \ell_{u}(n+m).  \notag
\end{equation}
\end{itemize}
\end{fact}

\begin{fact} \label{fact:2-5} 
 (1) Let $u\in C_{+, \log}$ be (log, exp)-convex. Then
 its $L$-function $\cl_{u}$ is also (log, exp)-convex and 
 for any $a>1$,
 \begin{equation} 
  \cl_{u}(r) \leq {ea \over \log a} u(ar), \qquad
  \forall r\geq 0.  \notag
 \end{equation}
 \par\noindent
 (2) Let $u\in C_{+, \log}$ be increasing and 
 (log, $x^{k}$)-convex. Then there exists a constant $C$, 
 independent of $k$, such that
 \begin{equation} 
  u(r) \leq C \cl_{u}(2^{k} r), \qquad \forall r\geq 0.
     \notag
 \end{equation}
 \end{fact}

\begin{fact} \label{fact:2-6}
Let $u\in C_{+, 1/2}$. Then its dual Legendre transform 
$u^{*}$ belongs to $C_{+, 1/2}$ and is an increasing 
(log, $x^{2}$)-convex function on $[0, \infty)$.
\end{fact}

\begin{fact} \label{fact:2-7}
If $u\in C_{+, 1/2}$ is (log, $x^{2}$)-convex, then the
Legendre transform $\,\ell_{u^{*}}$ of $u^{*}$ is given by
\begin{equation}
 \ell_{u^{*}} (t) = {e^{2t} \over \ell_{u}(t) t^{2t}},
 \qquad t\in [0, \infty).   \notag
\end{equation}
\end{fact}

\begin{fact} \label{fact:2-8}
Let $u\in C_{+, 1/2}$ be (log, $x^{2}$)-convex. If $u$
is increasing on the interval $[r_{0}, \infty)$, then we 
have $(u^{*})^{*}(r) = u(r)$ for all $r\geq r_{0}$. In 
particular, if $u$ is increasing on $[0, \infty)$, then
$(u^{*})^{*} = u$ on $[0, \infty)$.
\end{fact}

\begin{fact} \label{fact:2-9}
Let $u \in C_{+, 1/2}$ be (log, $x^{2}$)-convex. Then the
functions $u^{*}, \>\cl_{u^{*}}$, and $\cl_{u}^{\#}$ are
all equivalent.
\end{fact}

\smallskip
\begin{lemma} \label{lem:2-1}
The conditions (B1) and ($\wt{\text{B}}1$) are equivalent to 
near-(B2) and near-($\wt{\text{B}}2$), respectively. 
\end{lemma}

\begin{pf} 
It is enough to show the equivalence of (B1) and near-(B2). 
Put $u(r) = G_\alpha(r)$.  It is easy to see that 
$u \in C_{+,\log}$ and  $\alpha(n)/n! \leq \ell_u(n)$. 
By Fact \ref{fact:2-1}, $\ell_u(n)$ is log-concave. 
Since $\inf_{r>0} G_\alpha(r)/r^n = \ell_u(n)$, 
the condition (B1) is equivalent to that there exists 
a positive constant $C$ such that 
$\ell_u(n) \leq C^n \alpha(n)/n!$.  Hence $\{\a(n)/n!\}$ 
is equivalent to the log-concave sequence $\{\ell_u(n)\}$,
if (B1) holds. 

Conversely, suppose that there exists a positive sequence 
$\{\beta(n)\}$ and positive constants 
$K_1, K_2, c_1, c_2$ such that 
$\{\beta(n)/n!\}$ is log-concave. 
$$
K_1c_1^n{\beta(n) \over n!} \leq {\alpha(n) \over n!} 
\leq K_2c_2^n{\beta(n) \over n!}.
$$
Then we have
$$
K_1G_\beta(c_1r) \leq G_\alpha(r) \leq K_2G_\beta(c_2r).
$$
Therefore, we obtain
\begin{align}
\limsup_{n\to\infty} \left({n! \over \a(n)}
 \inf_{r>0} {G_{\a}(r) \over r^{n}}\right)^{1/n} 
& \leq \limsup_{n\to\infty} \left({n! \over K_1c_1^n\beta(n)} 
 \inf_{r>0} {K_2G_{\beta}(c_2r) \over r^{n}}\right)^{1/n} 
\notag \\
& = {K_2c_2 \over K_1c_1}
\limsup_{n\to\infty} \left({n! \over \beta(n)}
 \inf_{r>0} {G_{\beta}(r) \over r^{n}}\right)^{1/n} < \infty, 
\notag
\end{align}
since the condition (B2) for $\{\beta(n)\}$ implies 
(B1) for $\{\beta(n)\}$ (see \cite{cks}). 
\end{pf}

\smallskip
For more precise discussion,  We will need the 
following conditions on $u \in C_{+, 1/2}$: 
\begin{itemize}
\smallskip
\item[(U0)] $\inf_{r\geq 0} u(r) = 1$.
\smallskip
\item[(U1)] $u$ is increasing and $u(0)=1$.
\smallskip
\item[(U2)] $\lim_{r\to\infty} r^{-1} \log u(r) < \infty$.
\smallskip
\item[(U3)] $u$ is (log, $x^{2}$)-convex.
\end{itemize}

Recall that the A-conditions are needed in 
order to set up the Gel'fand triple $[\ce]_{\a} \subset 
(L^{2}) \subset [\ce]_{\a}^{*}$ and to make sure that
the renormalized exponential functions $:\!e^{\la\cdot,
\x\ra}\!:,\, \x\in\ce_{c}$, are test functions in
$[\ce]_{\a}$. Moreover, note that 
the B-conditions are used for the
characterization theorems \cite{akk2} \cite{cks}.
Keeping these in mind,
we consider the relationship between the 
U-conditions and AB-conditions in the rest of this section.

For a given $u\in C_{+,\log}$, define a sequence
$\{\a_u(n)\}$ by 
\begin{equation}\label{eq:2-4}
	\a_u(n)={1\over n!\ell_u(n)}.
\end{equation}

\begin{lemma} \label{lem:3-1}
If $u\in C_{+,\log}$ satisfies conditions $(U0)$ and 
$(U2)$, then $\alpha_u$ satisfies 
the condition (A1).
\end{lemma}

\begin{pf}
By the definition of Legendre transform, $\ell_u(0) = 1$.
Since $u$ satisfies condition (U2), there exist constants 
$c_{1}, c_{2}>0$ such that $u(r) \leq c_{1}e^{c_{2} r}$ 
for all $r\geq 0$. Therefore,
\begin{equation*} 
\ell_{u}(n) = \inf_{r>0} {u(r) \over r^{n}} \leq
   \inf_{r>0} {c_{1} e^{c_{2} r} \over r^{n}} = 
   c_{1} c_{2}^{n} \left({e\over n}\right)^{n} 
   \leq c_1 e {(c_2\sqrt{2}\,)^{n} \over n!}. 
\end{equation*}
by the Stirling formula $n!\leq e\sqrt{n}\big(n/e\big)^{n}$.
Therefore $\alpha_u(n) (c_2\sqrt{2}\,) \geq (c_1e)^{-1}$. 
\end{pf}

Further, we can show the condition (A2) by the following lemma. 

\begin{lemma} \label{lem:3-2}
Let $u\in C_{+, 1/2}$ satisfy condition (U3) and 
$\alpha_u$ be in (\ref{eq:2-4}). 
Then, $\alpha_u$ satisfies the condition (A2).  
In addition, $G_{\alpha_u} $ defined in Equations (\ref{eq:1-8}) 
and ${\cal L}_{u}^{\#}$ in  (\ref{eq:2-2}) 
are the same entire function, i.e. 
$G_{\alpha_u}(r) = {\cal L}_{u}^{\#}(r)$. 
\end{lemma} 

\begin{pf}
The equality is obvious. 
By the condition $u\in C_{+, 1/2}$ and by 
Fact \ref{fact:2-6} the dual transform $u^{*}$ belongs 
to $C_{+, 1/2}$ and is an increasing (log, $x^{2}$)-convex 
function. By Fact \ref{fact:2-0}, $u^{*}$ belongs to 
$C_{+, \log}$ and is  (log, exp)-convex. Therefore 
$\ell_{u^*}$ satisfies (2) in Fact \ref{fact:2-3}.
Then we see 
\begin{equation*} 
\lim_{n\to\infty} \left({1 \over \ell_{u}(n)(n!)^2}
\right)^{1/n}
= \lim_{n\to\infty} \left({\ell_{u^*}(n) n^{2n} \over (n!)^2 e^{2n}}
\right)^{1/n}
=  \lim_{n\to\infty} \ell_{u^{*}}(n)^{1/n} = 0 
\end{equation*}
by Fact \ref{fact:2-7} and the Stirling formula we have  
the condition (A2).
\end{pf}

\begin{lemma} \label{lem:5-3}
If $u\in C_{+, \,\log}$ satisfies condition (U3), then
the sequence $\{\a_{u}(n)\}$ satisfies condition 
near-(B2).
\end{lemma}

\noindent
{\em Remark.} Even if we assume that $u\in C_{+, 1/2}$
with condition (U3), we cannot conclude that $\{\a_{u}(n)\}$
satisfies condition (B2). On the other hand, if we assume
that $u\in C_{+, \log}$ is log-convex, then $\{\a_{u}(n)\}$
does satisfy condition (B2). For the proof, see Lemma 3.3
in \cite{akk4}.

\begin{pf}
We can apply Fact \ref{fact:2-4} (1) to see that
$\{\ell_{u}(n)n^{2n}\}$ is log-convex. However $\ell_{u}(n)=
(\a_{u}(n) n!)^{-1}$. Hence the sequence $\{(\a_{u}(n) 
n!)^{-1} n^{2n}\}$ is log-convex and so the sequence 
$\{\a_{u}(n)n!/n^{2n}\}$ is log-concave. 

Let $\l(n)=\a_{u}(n)(n!)^{2}/n^{2n}$. We have just shown 
that $\{\l(n)/n!\}$ is log-concave. On the other hand, it 
follows from the Stirling formula that $\{\a_{u}(n)\}$ and
$\{\l(n)\}$ are equivalent. Hence $\{\a(n)\}$ satisfies 
condition near-(B2).
\end{pf}

\begin{lemma} \label{lem:5-4}
Let $u\in C_{+, \log}$. Then the sequence $\{\a_{u}(n)\}$
satisfies condition ($\wt{\text{B}}2$).
\end{lemma}

\begin{pf}
By Fact \ref{fact:2-1}, $\{\ell_{u}(n)\}$ is 
log-concave. However $\ell_{u}(n)=(n!\a_{u}(n))^{-1}$ and so
the sequence $\{(n!\a_{u}(n))^{-1}\}$ is log-concave. This
means that the sequence $\{\a_{u}(n)\}$ satisfies condition 
($\wt{\text{B}}2$).
\end{pf}

Putting the above four lemmas together,
we get the next theorem for 
the Gel'fand triple $[\ce]_{u}\subset (L^{2}) \subset 
[\ce]_{u}^{*}$ associated with a growth function $u$.

\begin{theorem} \label{thm:5-1}
Suppose $u\in C_{+, 1/2}$ satisfies conditions (U0) (U2) 
(U3). Then the sequence $\a_{u}(n) = \big(\ell_{u}(n) 
n!\big)^{-1}, \> n\geq 0$, satisfies conditions (A1), (A2), 
near-(B2), and ($\wt{\text{B}}2$).
\end{theorem}

\smallskip
\section{Characterization theorems}
    \label{sec:3}

\smallskip

In the following we construct a Gel'fand triple 
$[\ce]_{u} \subset (L^{2}) \subset [\ce]_{u}^{*}$ 
associated with $u\in C_{+, 1/2}$ satisfying 
conditions (U0) (U2) (U3) and discuss characterization 
theorems  for test and generalized functions 
under the same condition. 

Note that condition (U0) is merely a normalization condition 
and is equivalent to $\ell_{u}(0)=1$, which guarantees the
condition $\alpha(0)=1$ in (A1).  Obviously, (U1) is 
stronger than condition (U0).  However we have 
the following lemma:

\begin{lemma}\label{lem:3}
For $u\in C_{+, 1/2}$, there exists  
a minimum point  $\underline{r}$ of $u$ and a maximum 
point $\bar{r}$ of $u$ on $[0,\underline{r}]$,
i.e. $u(\bar{r}) = 
\inf_{r \geq 0}u(r)$ and $u(\bar{r}) 
= \sup_{0 \leq r \leq \underline{r}}u(r)$. Define 
$$
v(r) = 
\begin{cases}
      u(\underline{r}) \text{\ for $0 \leq r \leq \underline{r}$}\\ 
      u(r) \text{\ for $\underline{r} \leq r$}.
\end{cases}
$$
Then $v$ belongs to $C_{+, 1/2}$ and equivalent to $u$;
$$
v(r) \leq u(r) \leq {u(\bar{r}) \over u(\underline{r})}v(r).
$$ 
Moreover 
$$
\ell_v(t) = \ell_u(t) \ \ \hbox{for } \ t \geq 0 \quad
\hbox{and } \quad \cl_v(r) = \cl_u(r) \ \ \hbox{for } \ r \geq 0 .
$$
If $u$ satisfies (U0) and (U3), then $v$ satisfies 
(U1)  and (U3) with $v(0) = u(\underline{r}) =1 $ and 
$u(\bar{r})=u(0)$.
\end{lemma}

\begin{pf}
Since $r^{-t}$ is decreasing on $[0,\ \underline{r}]$ 
for a fixed $t\geq 0$ and $u(r)$ takes the minimum value 
at $\underline{r}$, 
$$
\inf_{0 < r \leq \underline{r}}u(r)r^{-t} 
= u(\underline{r})\underline{r}^{-t} 
= v(\underline{r})\underline{r}^{-t} = 
\inf_{0 < r \leq \underline{r}}v(r)r^{-t}.
$$
This implies $\ell_v(t) = \ell_u(t)$ for $t \geq 0$. 
Other statements are obvious from above. 
\end{pf}

If $u$ satisfies (U0) and (U3), then 
by Fact \ref{fact:2-0}, this $v$ is $(\log, \exp)$-convex and 
hence  by Fact \ref{fact:2-5} $v$ is equivalent to $\cl_v$.   
This means that $u$ is equivalent to $\cl_u$ by the above lemma.

Now we will construct a Gel'fand triple 
$[\ce]_{u} \subset (L^{2}) \subset [\ce]_{u}^{*}$ 
associated with  a fixed function 
$u \in C_{+,1/2}$ satisfying 
conditions (U0) (U2) (U3).  
First we will relate the Gel'fand
triple $[\ce]_{u} \subset (L^{2}) \subset [\ce]_{u}^{*}$
to a CKS-space.
By (U3) and Fact \ref{fact:2-0}, we have a log-convex function 
$\{\ell_u(t)\}$. Define log-concave sequence $\alpha_u$ 
like \eqref{eq:2-4} by 
\begin{equation}\label{eq:alpha}
\alpha_u(n) = {1 \over n!\ell_u(n)}.
\end{equation}

Thus, by Lemmas \ref{lem:3-1} and \ref{lem:3-2} 
we can construct a Gel'fand triple 
$[\ce]_{u_\alpha} \subset (L^{2}) \subset [\ce]_{u_\alpha}^{*}$, 
which is denoted by 
$[\ce]_{u} \subset (L^{2}) \subset [\ce]_{u}^{*}$. 
Norms for $p > 0$ are also denoted as 
\begin{equation}\label{eq:3-7}
\|\vf\|_{p, u}= \|\vf\|_{p,\alpha_u}
= \left(\sum_{n=0}^{\infty} 
  {1\over \ell_{u}(n)}\,|f_{n}|_{p}^{2}\right)^{1/2}
\end{equation}
for $\vf=\sum_{n=0}^{\infty} \la :\!\cdot^{\otimes n}\!:, 
f_{n}\ra \in [\ce]_u$ and   
\begin{equation}\label{eq:3-8}
\|\F\|_{-p, (u)}= \|\F\|_{-p,1/\alpha_u}.
= \left(\sum_{n=0}^{\infty} \ell_{u}(n)
  (n!)^{2} |f_{n}|_{-p}^{2}\right)^{1/2}
\end{equation}
for $\F=\sum_{n=0}^{\infty} \la :\!\cdot^{\otimes n}\!:, 
f_{n}\ra \in [\ce]_u^*$. 
Subspaces are denoted as 
$$
[\ce_p]_u = [\ce_p]_\alpha \quad\hbox{and}\quad
[\ce_p]_u^* = [\ce_p]_\alpha^* \quad \hbox{for \ } \ p > 0. 
$$

Applying Theorem \ref{thm:1-1}, we can show the 
following theorem: 

\begin{theorem} \label{thm:3-1}
Suppose $u\in C_{+, 1/2}$ satisfies conditions (U0) (U2) 
(U3).   Then

\noindent
(i) \ There exist positive constants $c$ and $a$ 
such that for $\F \in [\ce]_{u}^{*}$
\begin{equation} \label{eq:gg}
 |(S\F)(\x)| \leq \|\F\|_{-p, (u)} \sqrt{c}\,
  u^{*}\big(a|\x|_{p}^{2}\big)^{1/2}, 
  \qquad \x\in\ce_{c}. 
\end{equation}

\noindent
(ii) \ A complex-valued function $F$ on $\ce_{c}$ is
the $S$-transform of a generalized function $\F \in
[\ce]_{u}^{*}$ if and only if it satisfies the conditions:
\begin{itemize}
\item[(1)] For any $\x, \y \in \ce_{c}$, the function
$F(z\x+\y)$ is an entire function of $z\in\spc$.
\item[(2)] There exist constants $K, a, p\geq 0$ such
that
\begin{equation} 
 |F(\x)| \leq K u^{*}\big(a|\x|_{p}^{2}\big)^{1/2},
  \qquad \x\in\ce_{c}.  \notag
\end{equation}
\end{itemize}
\noindent
(iii) \ In the above case, for any $q>p$ such that $ae^{2}
\|i_{q, p}\|_{HS}^{2}<1$, we have the inequality:
\begin{equation} \label{eq:3-13}
 \|\F\|_{-q, (u)} \leq K\left(1-ae^{2}\|i_{q, p}\|_{HS}^{2}
 \right)^{-1/2}. 
\end{equation}
\end{theorem}

\noindent
{\em Remark.} The growth condition (2) is equivalent
to the condition: There exist constants $K, p\geq 0$ such
that
\begin{equation}
 |F(\x)| \leq K u^{*}\big(|\x|_{p}^{2}\big)^{1/2},
  \qquad \x\in\ce_{c}.  \notag
\end{equation}

\begin{pf}
By Equation (\ref{eq:1-5}) and Lemma \ref{lem:3-2}, we have 
$$
 \|\!:\!e^{\la \cdot, \x\ra}\!:\!\|_{p, u}
= \cl_{u}^{\#}\big(|\x|_{p}^{2}\big)^{1/2} \quad p \geq 0. 
$$
Since $G_{\alpha_u}(r) = {\cal L}_{u}^{\#}(r)$. 
is equivalent to $u^*$ by Lemma \ref{lem:3-2}, 
$$
|S\F(\xi)| = |\la \F, :\!e^{\la \cdot, \x\ra}\!:\ra |
\leq \|\F\|_{-p,(u)} \cl_{u}^{\#}\big(|\x|_{p}^{2}\big)^{1/2}
\leq \|\F\|_{-p,(u)} \sqrt{c}u^*\big(a|\x|_{p}^{2}\big)^{1/2}
$$
with suitable $c, a > 0$. Thus we see (i).
By Lemma \ref{lem:5-3}, $\{\a(n)\}$ satisfies
condition near-$(B2)$.
By Fact \ref{fact:2-9}, (ii) follows from Theorem \ref{thm:1-1}.

  We can prove the estimation of the norm (\ref{eq:3-13}) 
with the same idea as in the proof of Theorem 8.2 in \cite{kuo96}. 
$\F$ and $F=S\F$ are expanded as 
$\F = \sum_{n=0}^{\infty} \la :\!\cdot^{\otimes n}\!:, f_{n}\ra$.
$\displaystyle  F(\x) = \sum_{n=0}^{\infty} \la f_{n}, 
  \x^{\otimes n}\ra$, respectively.
For any $\x_{1}, \ldots, \x_{n} \in \ce_{c}$,  applying 
the Cauchy formula to 
$\displaystyle F(z_1\xi_1+\cdots+z_n\xi_n)$, we have 
$$
 |\la f_{n}, \x_{1}\wh\otimes \cdots\wh\otimes 
    \x_{n}\ra| \leq {K \over n!} a^{n/2} n^{n}
    \left({u^{*}(r) \over r^{n}}\right)^{1/2} 
    |\x_{1}|_{p} \cdots |\x_{n}|_{p}.  
$$
Take the infimum over $r>0$ and use the definition of the 
Legendre transform.
Then, we conclude that
\begin{equation} \label{eq:3-14}
 |f_{n}|_{-q}^{2} \leq {K^{2} \over (n!)^{2}} a^{n} n^{2n}
  \ell_{u^{*}}(n) \|i_{q, p}\|_{HS}^{2n}. 
\end{equation}
Then by Equations (\ref{eq:3-8}) 
and (\ref{eq:3-14}),
$$
 \|\F\|_{-q, (u)}^{2} 
  = \sum_{n=0}^{\infty} \ell_{u}(n)(n!)^{2} 
     |f_{n}|_{-q}^{2}  \notag  \\
  \leq K^{2} \sum_{n=0}^{\infty} \ell_{u}(n)
     a^{n} n^{2n} \ell_{u^{*}}(n) 
       \|i_{q, p}\|_{HS}^{2n}.  
$$
However by Fact \ref{fact:2-7}, we have $\ell_{u^{*}}(n)=
\ell_{u}(n)^{-1} n^{-2n} e^{2n}$. Therefore,
\begin{equation}
 \|\F\|_{-q, (u)}^{2} \leq K^{2} \sum_{n=0}^{\infty}
  a^{n} e^{2n} \|i_{q, p}\|_{HS}^{2n}
 \leq K^2\left(1-ae^{2}
  \|i_{q, p}\|_{HS}^{2}\right)^{-1}  \notag
\end{equation}
by the assumption $ae^{2}\|i_{q, p}\|_{HS}^{2}<1$. 
(Of course, this estimation implies (ii), directly.)
\end{pf}

\begin{proposition} \label{prop:3-1}
(i)
For $\F=\sum_{n=0}^{\infty} \la 
:\!\cdot^{\otimes n}\!:, f_{n}\ra$ and $p\geq 0$, we define a new norm
\begin{equation} \label{eq:3-9}
 \|\F\|_{-p, u^{*}} = \left(\sum_{n=0}^{\infty}
  {1\over \ell_{u^{*}(n)}}\,|f_{n}|_{-p}^{2}
    \right)^{1/2}. 
\end{equation}
For any $p\geq 0$ and $q\geq p+{\log 2 \over 2\log
(\r^{-1})}$, we have
\begin{equation} \label{eq:3-10}
 e^{-1}\, \|\F\|_{-q, (u)} \leq \|\F\|_{-p, u^{*}}
 \leq \|\F\|_{-p, (u)}, \qquad 
   \forall \F \in [\ce_{p}]_{u}^{*}.
\end{equation}
In particular, if $S\F(\xi)$ satisfies the conditions (1) and (2) 
in Theorem \ref{thm:3-1}, we have
\begin{equation}
 \|\F\|_{-q, u^{*}} \leq K\left(1-ae^{2}\|i_{q, p}\|_{HS}^{2}
 \right)^{-1/2}.  \notag
\end{equation}

\noindent
(ii)
For a test function $\vf=\sum_{n=0}^{\infty} \la 
:\!\cdot^{\otimes n}\!:, f_{n}\ra$ and $p\geq 0$,  a new norm
corresponding to the norm $\|\cdot\|_{-p, u^{*}}$ in Equation 
(\ref{eq:3-9}) is given by
\begin{equation} \label{eq:3-b}
 \|\vf\|_{p, (u^{*})} = \left(\sum_{n=0}^{\infty}
  \ell_{u^{*}(n)} (n!)^{2} |f_{n}|_{p}^{2}
    \right)^{1/2}. 
\end{equation}
We have the corresponding inequalities, i.e., for any
$p\geq 0$ and $q\geq p+{\log 2\over 2\log (\r^{-1})}$,
\begin{equation}
\|\vf\|_{p, u} \leq \|\vf\|_{p, (u^{*})} \leq
 e\,\|\vf\|_{q, u}, \qquad 
  \forall \vf\in [\ce_{q}]_{u}.  \notag
\end{equation}
\end{proposition}

\begin{pf}
First we point out the following inequalities from page
357 in \cite{kuo96}
\begin{equation} \label{eq:3-11}
 e^{-1} 2^{-n/2} n! \leq \left({n\over e}\right)^{n}
  \leq n!.
\end{equation}
By Fact \ref{fact:2-7} and the second inequality in 
Equation (\ref{eq:3-11}),
\begin{equation}
 \|\F\|_{-p, u^{*}}^{2} = \sum_{n=0}^{\infty} \ell_{u}(n)
 \left({n\over e}\right)^{2n} |f_{n}|_{-p}^{2}
 \leq \sum_{n=0}^{\infty} \ell_{u}(n) (n!)^{2} 
  |f_{n}|_{-p}^{2} = \|\F\|_{-p, (u)}^{2}.  \notag
\end{equation}
This gives the second inequality in Equation 
(\ref{eq:3-10}). On the other hand, we can use the first 
inequality in Equation (\ref{eq:3-11}) to get
\begin{equation}
 \|\F\|_{-p, u^{*}}^{2} \geq e^{-2} \sum_{n=0}^{\infty}
 \ell_{u}(n) 2^{-n} (n!)^{2} |f_{n}|_{-p}^{2}.  \notag
\end{equation}
Note that $|f|_{-p}\geq \r^{p-q}|f|_{-q}$ for any $q\geq 
p$ and $f\in\ce_{p}'$. Therefore, 
\begin{equation}
 \|\F\|_{-p, u^{*}}^{2} \geq e^{-2} \sum_{n=0}^{\infty}
 \ell_{u}(n) (n!)^{2} \left(2^{-1} \r^{-2(q-p)}\right)^{n}
   |f_{n}|_{-q}^{2}.  \notag
\end{equation}
When $2\r^{2(q-p)}\leq 1$, i.e., $q\geq p+{\log 2 \over 
2\log (\r^{-1})}$, the above inequality yields that
\begin{equation}
 \|\F\|_{-p, u^{*}}^{2} \geq e^{-2} \sum_{n=0}^{\infty}
 \ell_{u}(n) (n!)^{2} |f_{n}|_{-q}^{2}
  = e^{-2}\|\F\|_{-q, (u)}^{2}.  \notag
\end{equation}
This proves the first inequality in Equation 
(\ref{eq:3-10}).
The assertions to test functions are proved similarly.
\end{pf}

Next we consider the characterization of test functions.
\begin{theorem} \label{thm:3-2}
Suppose $u\in C_{+, 1/2}$ satisfies conditions (U0) (U2) (U3). 
Then 

\noindent
(i) \ There exists a positive constant $a$ such that 
for $\vf \in [\ce]_{u}$,
$$
|S\vf(\xi)| \leq \|\vf\|_{u,p}\cl_u\big(|\xi|_p^2\big)
\leq \|\vf\|_{u,p}\sqrt{2e \over \log 2} u\big(a|\xi|_p^2\big).
$$

\noindent
(ii) \ A complex-valued function $F$ on $\ce_{c}$ is
the $S$-transform of a test function $\vf \in [\ce]_{u}$ 
if and only if it satisfies the conditions:
\begin{itemize}
\item[(1)] For any $\x, \y \in \ce_{c}$, the function
$F(z\x+\y)$ is an entire function of $z\in\spc$.
\item[(2)] There exist constants $K, a, p\geq 0$ such
that
\end{itemize}
\begin{equation} \label{eq:3-18}
 |F(\x)| \leq K u\big(a|\x|_{-p}^{2}\big)^{1/2},
  \qquad \x\in\ce_{c}. 
\end{equation}

\noindent
(iii) \ Let $q\in [0, p)$ be a number such that 
$ae^{2}\|i_{p, q}\|_{HS}^{2} < 1$. Then $F$ is the 
$S$-transform of $\vf\in [\ce_{q}]_{u}$ and 
\begin{equation} \label{eq:3-a}
 \|\vf\|_{q, u} \leq K\left(1-ae^{2}\|i_{p, q}\|_{HS}^{2}
 \right)^{-1/2}. 
\end{equation}
\end{theorem}

\noindent
{\em Remark.} The growth condition (2) is equivalent to 
the condition: For any $p\geq 0$ there exists a constant
$K\geq 0$ such that
\begin{equation}
 |F(\x)| \leq K u\big(|\x|_{-p}^{2}\big)^{1/2},
  \qquad \x\in\ce_{c}.  \notag
\end{equation}

\begin{pf}
Let $v$ be as in Lemma \ref{lem:3}. 
It is easy to see that
$G_{1/\alpha_u} (r) = \cl_u(r) = \cl_v(r)$ and hence
$$
 \|\!:\!e^{\la \cdot, \x\ra}\!:\!\|_{p,(u)}
= \cl_{u}\big(|\x|_{-p}^{2}\big)^{1/2} \quad p \geq 0
$$
by Equation (\ref{eq:1-7}). Therefore
$$
|S\vf(\xi)| = |\la \F, :\!e^{\la \cdot, \x\ra}\!:\ra |
\leq \|\vf\|_{p,u} \cl_{u}\big(|\x|_{-p}^{2}\big)^{1/2}.
$$
By using Fact \ref{fact:2-5} (1) for $v$ with $a=2$, we have
$$
\cl_u(r) = \cl_v(r) \leq {2e\over \log 2}v(2r) 
\leq {2e\over \log 2}u(2r).
$$ 
Thus we see (i).

By using Fact \ref{fact:2-5} (2) for $v$ with $k=2$, we have
$$
u(r) \leq u(0)v(r) 
\leq Cu(0)\cl_v(2^2r) 
= Cu(0)\cl_u(2^2r) .
$$
We have already seen $\cl_u(r) \leq {2e\over \log 2}u(2r)$ 
in the proof of (i). Thus $u$ and 
$\cl_u = G_{1/\alpha_u}$ are equivalent. 
Due to Lemma \ref{lem:5-4},
$\{\a_u(n)\}$ satisfies the condition
($\wt{\text{B}}2$). 
By Theorem \ref{thm:1-2}, we can prove (ii).

   The proof of (iii) is similar to that of
Theorem \ref{thm:3-1} (iii).  
The key of the proof is 
the estimation
$$
|f_n|_p \leq \inf_{r>0} {u(an^2r^2)^{1/2} \over r^n}
= a^{n/2}n^n \ell_u(n)^{1/2}
\leq (a^{1/2}e)^n n! \ell_u(n)^{1/2}
$$
for $\displaystyle  F(\x) = \sum_{n=0}^{\infty} \la f_{n},
\x^{\otimes n}\ra$.  This implies (iii).
\end{pf}

The following examples
can be applicable to 
Theorems \ref{thm:3-1} and \ref{thm:3-2}.
\begin{example}
Consider 
\begin{equation*}
	u(r)=u^{*}(r)=e^r.
\end{equation*}
Then it is obvious to check 
that consditions (U0) (U2) (U3) are satisfied.
\end{example}

\begin{example}
For $0\leq \b<1$, let $u$ be the function defined by
\begin{equation}
  u(r)= \exp\left[(1+\b)r^{1\over 1+\b}\right].  \notag
\end{equation}
It is easy to check that $u$ belongs to $C_{+, 1/2}$ and 
satisfies conditions (U0) (U2) (U3). By Example 4.3 in 
\cite{akk3}, the dual Legendre transform $u^{*}$ of $u$ 
is given by
\begin{equation}
 u^{*}(r) = \exp\left[(1-\b)r^{1\over 1-\b}\right]. \notag
\end{equation}
Hence Theorems \ref{thm:3-1} and \ref{thm:3-2}
can be applied to the Gel'fand triple 
$[\ce]_{u}\subset (L^{2}) \subset [\ce]_{u}^{*}$ for the 
pair of functions $u^{*}$ and $u$.
\end{example}

\begin{example}
Consider the function $v(r) = \exp\big[e^{r}-1\big]$. 
Obviously, $v\in C_{+, 1/2}$. Let $u=v^{*}$ be the dual
Legendre transform of $v$. Then $u(0)=\sup_{s\geq 0}
v(s)^{-1}=1$ and by Fact \ref{fact:2-6} $u$ belongs to
$C_{+, 1/2}$ and is an increasing (log, $x^{2}$)-convex 
function on $[0, \infty)$. Hence $u\in C_{+, 1/2}$ 
satisfies conditions (U1) and (U3). It is shown in Example
4.4 in \cite{akk3} that $u$ is equivalent to the function
\begin{equation} \label{eq:3-20}
 w(r) = \exp\left[2\sqrt{r\log\sqrt{r}}\,\right]. 
\end{equation}
Obviously, $w$ satisfies condition (U2) and so $u$ also
satisfies condition (U2). On the other hand, we have 
$u^{*}=(v^{*})^{*}=v$ by Fact \ref{fact:2-8}. Hence 
Theorems \ref{thm:3-1} and \ref{thm:3-2}  
can be applied to the Gel'fand triple $[\ce]_{u}\subset 
(L^{2}) \subset [\ce]_{u}^{*}$ for the following pair of
functions:
\begin{equation}
 u^{*}(r) = \exp\big[e^{r}-1\big], \quad
 u(r) = (u^{*})^{*}.  \notag
\end{equation}
Observe that there is no exact form for the function
$u$. Thus for Theorems \ref{thm:3-2} 
we should use the equivalent function $w$ in Equation
(\ref{eq:3-20}) as the growth function.

In general, let $\exp_{k}(r)=\exp(\exp(\cdots(\exp(r))))$
be the $k$-th iterated exponential function and
consider the function
\begin{equation}
 v_{k}(r) = {\exp_{k}(r) \over \exp_{k}(0)}.  \notag
\end{equation}
The dual Legendre transform $u_{k}=v_{k}^{*}$ belongs to
$C_{+, 1/2}$ and satisfies conditions (U1) (U2) (U3). The
function $u_{k}$ is equivalent to the function $w_{k}$
given in Equation (\ref{eq:1-11}), i.e.,
\begin{equation} \label{eq:3-21}
 w_{k}(r) = \exp\left[2\sqrt{r\log_{k-1}\sqrt{r}}\,\right].
\end{equation}
We have $u_{k}^{*}=v_{k}$ and Theorems \ref{thm:3-1},
and \ref{thm:3-2}                 
can be applied to the 
Gel'fand triple $[\ce]_{u_{k}}\subset (L^{2}) \subset 
[\ce]_{u_{k}}^{*}$ for the following pair of functions:
\begin{equation}
 u_{k}^{*}(r) = {\exp_{k}(r) \over \exp_{k}(0)}, \quad
 u_{k}(r) = (u_{k}^{*})^{*}.  \notag
\end{equation}
Again there is no exact form for the function $u_{k}$.
Thus for Theorems \ref{thm:3-2}  
we should use the equivalent function $w_{k}$ in Equation
(\ref{eq:3-21}) as the growth function.
\end{example}

\smallskip
\section{Intrinsic topology and Hida measures}
   \label{sec:4}

In the space $[\ce]_{u}$ of test functions there are two
families of norms, namely, $\{\|\cdot\|_{p, u};\,p\geq 0\}$ 
defined in Equation (\ref{eq:3-7}) and $\{\|\cdot\|_{p, 
(u^{*})};\,p\geq 0\}$ defined in Equation (\ref{eq:3-b}). 
As we pointed out in the Remark of Lemma \ref{prop:3-1},
these two families are equivalent. Observe that both
$\|\vf\|_{p, u}$ and $\|\vf\|_{p, (u^{*})}$ are defined
in terms of the Wiener-It\^o expansion of $\vf$.

In this section we will introduce another equivalent 
family of norms on $[\ce]_{u}$, i.e., 
$\{\|\cdot\|_{\ca_{p, u}};\,p\geq 0\}$. This family of 
norms is intrinsic in the sense that $\|\vf\|_{\ca_{p, u}}$ 
is defined directly in terms of the analyticity and growth 
condition of $\vf$.

First the analyticity, each test function $\vf$ in 
$[\ce]_{u}$ has a unique analytic extension (see \S 6.3 
of the book \cite{kuo96}) given by
\begin{equation} \label{eq:4-1}
 \vf(x) = \lla :\!e^{\la\cdot, x\ra}\!:, \T\vf\rra,
 \qquad x\in \ce_{c}', 
\end{equation}
where $\T$ is the unique linear operator taking 
$e^{\la\cdot, \x\ra}$ into $:\!e^{\la\cdot, \x\ra}\!:$ for 
all $\x\in \ce_{c}$. This operator turns out to be the 
same as $\cg_{i, 1}$ defined in Equation (\ref{eq:fg})
in section \ref{sec:5}. By Theorem \ref{thm:fg} this 
operator is continuous from $[\ce]_{u}$ into itself.

Now, let $p\geq 0$ be any fixed number. Choose $p_{1}>p$
such that $2\r^{2(p_{1}-p)}\leq 1$. Then use Equations
(\ref{eq:4-1}) and (\ref{eq:3-18}) to get
\begin{equation}
 |\vf(x)| \leq \|\T\vf\|_{p_{1}, u}\,\|\!:\!e^{\la\cdot,
  x\ra}\!:\!\|_{-p_{1}, (u)} 
  \leq \|\T\vf\|_{p_{1}, u}\, \sqrt{2e \over \log 2}\,
   u\big(2|x|_{-p_{1}}^{2}\big)^{1/2}.  \notag
\end{equation}
Note that $2|x|_{-p_{1}}^{2}\leq 2\r^{2(p_{1}-p)}
|x|_{-p}^{2} \leq |x|_{-p}^{2}$ by the above choice of 
$p_{1}$. Since $u$ is an increasing function, we see that
\begin{equation}
  |\vf(x)| \leq \|\T\vf\|_{p_{1}, u}\, \sqrt{2e \over 
  \log 2}\,u\big(|x|_{-p}^{2}\big)^{1/2}.  \notag
\end{equation}

However $\T$ is a continuous linear operator from $[\ce]_{u}$
into itself. Hence there exist positive constants $q$ 
and $K_{p, q}$ such that $\|\T\vf\|_{p_{1}, u}\leq 
K_{p, q}\|\vf\|_{q, u}$. Therefore,
\begin{equation} \label{eq:4-2}
  |\vf(x)| \leq C_{p, q} \|\vf\|_{q, u}\,u\big(|x|_{-p}^{2}
  \big)^{1/2},  \qquad x \in \ce_{p, c}',  
\end{equation}
where $C_{p, q}=K_{p, q}\sqrt{2e/\log 2}$. This is the
growth condition for test functions.

Being motivated by Equation (\ref{eq:4-2}), we define 
\begin{equation} \label{eq:4-3}
 \|\vf\|_{\ca_{p, u}} = \sup_{x\in\ce_{p, c}'} |\vf(x)|
 \,u\big(|x|_{-p}^{2}\big)^{-1/2}. 
\end{equation}
Obviously, $\|\cdot\|_{\ca_{p, u}}$ is a norm on 
$[\ce]_{u}$ for each $p\geq 0$. 
The next theorem generalizes the
results of Kuo \cite{kuo96} and Lee \cite{lee}.

\begin{theorem} \label{thm:4-1}
Suppose $u\in C_{+, 1/2}$ satisfies conditions (U0) (U2) 
(U3). Then the families of norms $\{\|\cdot\|_{\ca_{p, u}};
\,p\geq 0\}$ and $\{\|\cdot\|_{p, u};\,p\geq 0\}$ are 
equivalent, i.e., they generate the same topology on 
$[\ce]_{u}$.
\end{theorem}

\noindent
{\em Remark.}
(1) This theorem has been announced in \cite{akk5}, but 
(U1) condition was assumed there instead of (U0).
This theorem can be interpreted in another 
way which gives an alternative construction of test 
functions. \\
(2) For $p\geq 0$, let 
$\ca_{p, u}$ consist of all functions $\vf$ on $\ce_{c}'$ 
satisfying the conditions:
\begin{itemize}
\item[(a)] $\vf$ is an analytic function on $\ce_{p, c}'$.
\item[(b)] There exists a constant $C\geq 0$ such that
\begin{equation}
 |\vf(x)| \leq C u\big(|x|_{-p}^{2}\big)^{1/2}, \qquad
  \forall x \in \ce_{p, c}'.  \notag
\end{equation}
\end{itemize}
(3)
Suppose that $v$ is equivalent to $u$.  Then, 
it is obvious to see that the family of norms 
$\{\|\cdot\|_{\ca_{p,v}};\ p\geq 0\}$ and 
$\{\|\cdot\|_{\ca_{p,u}};\ p\geq 0\}$
are equivalent.

\medskip
For each $\vf\in \ca_{p, u}$, define $\|\vf\|_{\ca_{p, u}}$
by Equation (\ref{eq:4-3}). Then $\ca_{p, u}$ is a Banach 
space with norm $\|\cdot\|_{\ca_{p, u}}$. Let $\ca_{u}$ be 
the projective limit of $\{\ca_{p, u};\,p\geq 0\}$. We can
use the above theorem to conclude that $\ca_{u}=[\ce]_{u}$
as vector spaces with the same topology. Here the equality
$\ca_{u}=[\ce]_{u}$ requires the analytic extension of a
test function $\vf\in [\ce]_{u}$ in Equation 
(\ref{eq:4-1}).

\begin{pf}
Let $p\geq 0$ be any given number. We have already shown
that there exist constants $q>p$ and $C_{p, q}\geq 0$ 
such that Equation (\ref{eq:4-2}) holds. It follows that
\begin{equation}
 \|\vf\|_{\ca_{p, u}} = \sup_{x\in\ce_{p, c}'} |\vf(x)|
 \,u\big(|x|_{-p}^{2}\big)^{-1/2}
  \leq C_{p, q} \|\vf\|_{q, u}.  \notag
\end{equation}
Hence for any $p\geq 0$, there exist constants $q>p$ and 
$C_{p, q}\geq 0$ such that
\begin{equation} \label{eq:4-4}
  \|\vf\|_{\ca_{p, u}} \leq C_{p, q} \|\vf\|_{q, u},
  \qquad \forall \vf\in [\ce]_{u}.
\end{equation}

To show the converse, first note that by condition (U2)
there exist constants $c_{1}, c_{2}>0$ such that 
$u(r)\leq c_{1}e^{c_{2}r}, \> r\geq 0$. Next note that
by Fernique's theorem (see \cite{fer} \cite{kuo75} or
page 328 in \cite{kuo96}) we have
\begin{equation}
 \int_{\ce'} e^{2c_{2}|x|_{-\l}^{2}}\,d\m(x) < \infty
 \qquad \text{for all large~} \l.  \notag
\end{equation}

Now, let $p\geq 0$ be any given number. Choose $q>p$ 
large enough such that
\begin{equation} \label{eq:4-5}
 4e^{2}\|i_{q, p}\|_{HS}^{2}<1, \qquad
 \int_{\ce'} e^{2c_{2}|x|_{-q}^{2}}\,d\m(x) < \infty.
\end{equation}
With this choice of $q$ we will show below that
\begin{equation} \label{eq:4-6}
 \|\vf\|_{p, u} \leq L_{p, q} \|\vf\|_{\ca_{q, u}},
 \qquad \forall \vf\in [\ce]_{u},
\end{equation}
where $L_{p, q}$ is the constant given by
\begin{equation} \label{eq:4-a}
 L_{p, q} = \sqrt{c_{1}}\left(1-4e^{2}
  \|i_{q, p}\|_{HS}^{2}\right)^{-1/2} \int_{\ce'} 
  e^{2c_{2}|x|_{-q}^{2}}\,d\m(x). 
\end{equation}
Observe that the theorem follows from Equations 
(\ref{eq:4-4}) and (\ref{eq:4-6}).

To prove Equation (\ref{eq:4-6}), let $\vf\in [\ce]_{u}$ 
and $F=S\vf$. Then $F$ can be written as an integral 
(see page 36 in \cite{kuo96})
\begin{equation}
 F(\x) = \int_{\ce'} \vf(x+\x)\,d\m(x), \qquad
  \x\in\ce_{c}.  \notag
\end{equation}
Hence for the above choice of $q$, we have
\begin{align}
 |F(\x)| & \leq \int_{\ce'} |\vf(x+\x)|\,d\m(x)
      \notag  \\
  & \leq \int_{\ce'} \left(|\vf(x+\x)|\,
     u\big(|x+\x|_{-q}^{2}\big)^{-1/2}\right)
    u\big(|x+\x|_{-q}^{2}\big)^{1/2}\,d\m(x)
      \notag  \\
  & \leq \|\vf\|_{\ca_{q, u}} \int_{\ce'}
    u\big(|x+\x|_{-q}^{2}\big)^{1/2}\,d\m(x).
      \notag 
\end{align}
However by condition (U0), $u\geq 1$  
on $[0, \infty)$. Hence $u(r)^{1/2}\leq
u(r)$ for all $r\geq 0$. Therefore,
\begin{equation} \label{eq:4-7}
 |F(\x)| \leq \|\vf\|_{\ca_{q, u}} \int_{\ce'}
    u\big(|x+\x|_{-q}^{2}\big)\,d\m(x).
\end{equation}

By condition (U3), $u$ is (log, $x^{2}$)-convex. Thus
in particular, we have
\begin{equation}
 u\left(\big({\textstyle{1\over 2}}r_{1}+{\textstyle{1\over    
     2}}r_{2}\big)^{2}\right)
 \leq u\big(r_{1}^{2}\big)^{1/2}\,
   u\big(r_{2}^{2}\big)^{1/2},
 \qquad \forall r_{1}, r_{2} \geq 0. \notag
\end{equation}
Put $r_{1}=2|x|_{-q}$ and $r_{2}=2|\x|_{-q}$ to get
\begin{align}
 u\big(|x+\x|_{-q}^{2}\big) 
 & \leq u\left(\big({\textstyle{1\over 2}} 2|x|_{-q}
   +{\textstyle{1\over 2}} 2|\x|_{-q}\big)^{2}
       \right)      \notag \\
 & \leq u\big(4|x|_{-q}^{2}\big)^{1/2}\,
    u\big(4|\x|_{-q}^{2}\big)^{1/2}.   \notag
\end{align}
Then integrate over $\ce'$ to obtain the inequality:
\begin{equation} \label{eq:4-8}
 \int_{\ce'} u\big(|x+\x|_{-q}^{2}\big)\,d\m(x)
 \leq u\big(4|\x|_{-q}^{2}\big)^{1/2} \int_{\ce'}
  u\big(4|x|_{-q}^{2}\big)^{1/2}\,d\m(x).  
\end{equation}
Put Equation (\ref{eq:4-8}) into Equation (\ref{eq:4-7})
to get
\begin{equation} \label{eq:4-9}
 |F(\x)| \leq \|\vf\|_{\ca_{q, u}} 
  u\big(4|\x|_{-q}^{2}\big)^{1/2} \int_{\ce'} 
  u\big(4|x|_{-q}^{2}\big)^{1/2}\,d\m(x). 
\end{equation}

Now, by the inequality $u(r)\leq c_{1}e^{c_{2}r}$, we have
\begin{equation} \label{eq:4-10}
 \int_{\ce'} u\big(4|x|_{-q}^{2}\big)^{1/2}\,d\m(x)
   \leq  \sqrt{c_{1}} \int_{\ce'} e^{2c_{2}|x|_{-q}^{2}}
    \,d\m(x), 
\end{equation}
which is finite by the choice of $q$ in Equation 
(\ref{eq:4-5}).

By Equations (\ref{eq:4-9}) and (\ref{eq:4-10}),
we see that
\begin{equation}
  |F(\x)| \leq \|\vf\|_{\ca_{q, u}} \sqrt{c_{1}}
  \left(\int_{\ce'} e^{2c_{2}|x|_{-q}^{2}}\,d\m(x)\right)
    u\big(4|\x|_{-q}^{2}\big)^{1/2}, 
    \qquad \x\in\ce_{c}. \notag
\end{equation}
With this inequality and the choice of $q$ in Equation
(\ref{eq:4-5}) we can apply Theorem \ref{thm:3-2} to
show that for any $\vf\in [\ce]_{u}$,
\begin{equation}
 \|\vf\|_{p, u} \leq L_{p, q} \|\vf\|_{\ca_{q, u}}, \notag
\end{equation}
where $L_{p, q}$ is given by Equation(\ref{eq:4-a}). Thus
the inequality in Equation (\ref{eq:4-6}) holds and so the 
proof is completed.
\end{pf}

Next, we consider the characterization of Hida measures.
However first we need to prepare two lemmas.

\begin{lemma} \label{lem:4-1}
Suppose $u\in C_{+, \log}$ is (log, $x^{k}$)-convex. Then
\begin{equation}
 \cl_{u}(r)^{2} \leq \ell_{u}(0) \cl_{u}\big(2^{k+1}r\big),
  \qquad \forall r\in [0, \infty).
\end{equation}
\end{lemma}

\noindent
{\em Remark.} Note that $\cl_{u}(r)\geq \ell_{u}(0)$ for
all $r\geq 0$. Hence we have inequalities
\begin{equation}
 \ell_{u}(0) \cl_{u}(r) \leq \cl_{u}(r)^{2} \leq \ell_{u}(0) 
 \cl_{u}\big(2^{k+1}r\big), \qquad \forall r\in [0, \infty).
 \notag
\end{equation}
Thus $\cl_{u}$ and $\cl_{u}^{2}$ are equivalent for any 
(log, $x^{k}$)-convex function $u\in C_{+, \log}$. 
It follows that $u$ and
$u^{2}$ are equivalent for such a function $u$. 

\begin{pf}
Apply Fact \ref{fact:2-4} (2) to get
\begin{align}
 \cl_{u}(r)^{2} & = \sum_{j=0}^{\infty} \sum_{m=0}^{\infty}
     \ell_{u}(j) \ell_{u}(m) r^{j+m}   \notag  \\
  & \leq \ell_{u}(0) \sum_{j=0}^{\infty} \sum_{m=0}^{\infty}
    2^{k(j+m)} \ell_{u}(j+m) r^{j+m}  \notag  \\
  & = \ell_{u}(0) \sum_{j=0}^{\infty} \sum_{n=j}^{\infty}
       2^{kn} \ell_{u}(n) r^{n}.  \notag 
\end{align}
Then change the order of summation and use the inequality
$n+1 \leq 2^{n}$ to get
\begin{align}
 \cl_{u}(r)^{2} & \leq \ell_{u}(0) \sum_{n=0}^{\infty}
    (n+1) 2^{kn} \ell_{u}(n) r^{n}  \notag  \\
  & \leq \ell_{u}(0) \sum_{n=0}^{\infty}
     2^{(k+1)n} \ell_{u}(n) r^{n}  \notag  \\
  & = \ell_{u}(0) \cl_{u}\big(2^{k+1}r\big). \notag
\end{align}
\end{pf}

\begin{lemma} \label{lem:4-2}
Suppose $u\in C_{+, \log}$ is increasing and 
(log, $x^{k}$)-convex. Then for any $a>1$, we have
\begin{equation} \label{eq:4-11}
 \cl_{u}(r) \leq \sqrt{\ell_{u}(0){ea \over \log a}}\>
   u\big(a2^{k+1}r\big)^{1/2}.  
\end{equation}
\end{lemma}

\begin{pf}
Recall a fact mentioned in the beginning of section 
\ref{sec:2} that if $u$ is increasing and 
(log, $x^{k}$)-convex, then $u$ is (log, exp)-convex.
Hence this lemma follows from Lemma \ref{lem:4-1} and 
Fact \ref{fact:2-5} (1).
\end{pf}

A measure $\n$ on $\ce'$ is called a {\em Hida measure} 
associated with $u$ if $[\ce]_{u}\subset L^{1}(\n)$ and
the linear functional $\vf\mapsto \int_{\ce'} 
\vf(x)\,d\n(x)$ is continuous on $[\ce]_{u}$. In this case,
$\n$ induces a generalized function, denoted by $\wt\n$,
in $[\ce]_{u}^{*}$ such that
\begin{equation} \label{eq:4-12}
 \lla \wt\n, \vf\rra = \int_{\ce'} \vf(x)\,d\n(x),
  \qquad \vf \in [\ce]_{u}.
\end{equation}

The next theorem generalizes the
results of Kuo \cite{kuo96} and Lee \cite{lee}.
\begin{theorem} \label{thm:4-2}
Suppose $u\in C_{+, 1/2}$ satisfies conditions (U0) (U2) 
(U3). Then a measure $\n$ on $\ce'$ is a Hida measure with
$\wt\n\in [\ce]_{u}^{*}$ if and only if $\n$ is supported
by $\ce_{p}'$ for some $p\geq 0$ and 
\begin{equation} \label{eq:4-13}
 \int_{\ce_{p}'} u\big(|x|_{-p}^{2}\big)^{1/2}\,d\n(x)
  < \infty. 
\end{equation}
\end{theorem}

\noindent
{\em Remark.} 
This theorem has also been announced in \cite{akk5}, but
the conditions (U1) (U2) (U3) were assumed.

\begin{pf}
To prove the sufficiency, suppose $\n$ is supported by 
$\ce_{p}'$ for some $p\geq 0$ and Equation (\ref{eq:4-13}) 
holds. Then for any $\vf\in [\ce]_{u}$,
\begin{align}
 \int_{\ce'} |\vf(x)|\,d\n(x)
  & = \int_{\ce_{p}'} |\vf(x)|\,d\n(x)  \notag  \\
  & = \int_{\ce_{p}'} \left(|\vf(x)| u\big(|x|_{-p}^{2}
      \big)^{-1/2}\right) u\big(|x|_{-p}^{2}\big)^{1/2}
       \,d\n(x)  \notag  \\
  & \leq \|\vf\|_{\ca_{p, u}} \int_{\ce_{p}'}
   u\big(|x|_{-p}^{2}\big)^{1/2}\,d\n(x). \label{eq:4-14}
\end{align}

By Theorem \ref{thm:4-1}, $\{\|\cdot\|_{\ca_{p, u}};
\,p\geq 0\}$ and $\{\|\cdot\|_{p, u};\,p\geq 0\}$ are 
equivalent. Hence Equation (\ref{eq:4-14}) implies that
$[\ce]_{u}\subset L^{1}(\n)$ and the linear functional
\begin{equation} 
  \vf\longmapsto \int_{\ce'} \vf(x)\,d\n(x), \qquad
  \vf\in [\ce]_{u}, \notag
\end{equation}
is continuous on $[\ce]_{u}$. Thus $\n$ is a Hida measure 
with $\wt\n$ in $[\ce]_{u}^{*}$.

To prove the necessity, suppose $\n$ is a Hida measure 
inducing a generalized function $\wt\n \in [\ce]_{u}^{*}$. 
Then for all $\vf \in [\ce]_{u}$,
\begin{equation} \label{eq:4-15}
 \lla \wt\n, \vf \rra = \int_{\ce'} \vf(x)\,d\n(x).
\end{equation}
Since $\{\|\cdot\|_{\ca_{p, u}}; \,p\geq 0\}$ and 
$\{\|\cdot\|_{p, u};\,p\geq 0\}$ are equivalent, the
linear functional $\vf\mapsto \lla \wt\n, \vf \rra$ is
continuous with respect to $\{\|\cdot\|_{\ca_{p, u}}; 
\,p\geq 0\}$. Hence there exist constants $K, q\geq 0$
such that for all $\vf \in [\ce]_{u}$,
\begin{equation} \label{eq:4-16}
 \big| \lla \wt\n, \vf \rra\big| \leq 
   K\|\vf\|_{\ca_{q, u}}.
\end{equation}
Note that by continuity, Equations (\ref{eq:4-15}) and 
(\ref{eq:4-16}) also hold for all $\vf\in \ca_{q, u}$,
which is defined in the Remark of Theorem \ref{thm:4-1}.

Now, with this $q$, we define a function $\g$ on 
$\ce_{q, c}'$ by
\begin{equation}
 \g (x) = \cl_{u}\big(2^{-4} \la x, x\ra_{-q}\big),
 \qquad x\in \ce_{q, c}',  \notag
\end{equation}
where $\la \cdot, \cdot\ra_{-q}$ is the bilinear pairing
on $\ce_{q, c}'$. Obviously, $\g$ is analytic on 
$\ce_{q, c}'$. On the other hand, apply Lemma \ref{lem:4-2}
with $a=k=2$ to get
\begin{equation}
 |\g(x)| \leq \cl_{u}\big(2^{-4}|x|_{-q}^{2}\big) \leq 
  \sqrt{2e \over \log 2}\> u\big(|x|_{-q}^{2}\big)^{1/2}, 
   \qquad \forall x\in\ce_{q, c}'.  \notag
\end{equation}
This shows that $\g\in \ca_{q, u}$ and we have
\begin{equation}
 \|\g\|_{\ca_{q, u}} \leq \sqrt{2e \over \log 2}. \notag
\end{equation}
Then apply Equation (\ref{eq:4-16}) to the function $\g$,
\begin{equation}
  \big|\lla \wt\n, \g \rra\big| \leq K\|\g\|_{\ca_{q, u}} 
  \leq K \sqrt{2e \over \log 2}. \notag
\end{equation}
Therefore, from Equation (\ref{eq:4-15}) with $\vf=\g$
we conclude that
\begin{equation} \label{eq:4-17}
 \left|\int_{\ce'} \g(x)\,d\n(x)\right| \leq
   K \sqrt{2e \over \log 2}. 
\end{equation}

Note that $\g(x)=\cl_{u}\big(2^{-4}|x|_{-q}^{2}\big)$ for
$x\in\ce'$. Hence Equation (\ref{eq:4-17}) implies that
\begin{equation}
 \int_{\ce'} \cl_{u}\big(2^{-4}|x|_{-q}^{2}\big)\,d\n(x)
  < \infty.  \notag
\end{equation}
However $u(r) \leq C \cl_{u}(4r)$ from Fact \ref{fact:2-5} 
(2) with $k=2$. Therefore,
\begin{equation}
  \int_{\ce'} u\big(2^{-6}|x|_{-q}^{2}\big)\,d\n(x)
  < \infty.  \notag
\end{equation}
Now, choose $p>q$ large enough such that $\r^{2(p-q)}
\leq 2^{-6}$. Then $|x|_{-p}^{2} \leq 2^{-6}|x|_{-q}^{2}$.
Recall that $u$ is increasing. Hence
\begin{equation}
  \int_{\ce'} u\big(|x|_{-p}^{2}\big)\,d\n(x)
  < \infty.  \notag
\end{equation}
Note that $u(r)\geq 1$ and so $u(r)^{1/2} \leq u(r)$.
Thus we conclude that
\begin{equation}
 \int_{\ce'} u\big(|x|_{-p}^{2}\big)^{1/2}\,d\n(x)
  < \infty.  \notag
\end{equation}
This inequality implies that $\n$ is supported by 
$\ce_{p}'$ and Equation (\ref{eq:4-13}) holds.
\end{pf}

\smallskip
Before closing this section, 
let us explain the relationship with 
\cite{ghor}.  The basic equalities are 
\begin{equation*}
u(r) = e^{2\theta(\sqrt{r})}, \quad 
u^*(r) = e^{2\theta^*(\sqrt{r})}
\end{equation*}
where $\t^{*}(s)=\sup_{t>0}\{st-\t(t)\}$
is adopted in \cite{ghor}.
In the following table
we give the correspondence
between our $U$-conditions and $\t$-conditions.

\begin{table}[h]
\begin{center}
	\begin{tabular}{|l|c|c|}
	\noalign{\hrule height0.8pt}
	\hfil  \ & $u$ &
	$\theta$ \\
	\hline
	$(U0)$ & $\displaystyle \inf_{r\geq 0} u(r) = 1$ 
	& $\displaystyle \inf_{r\geq 0} \theta(r) = 0$ \\
	$(U1)$ & $u$ is increasing and $u(0)=1$
	& $\theta$ is increasing and 
	$\theta(0) = 0$\\
	$(U2)$ & $\displaystyle \lim_{r\to\infty} \frac{\log u(r)}{r} < \infty$
	& $\displaystyle \lim_{r\to\infty} \frac{\theta(r)}{r^2} < \infty$\\
	$(U3)$ & $u$ is (log, $x^{2}$)-convex 
	& $\theta$ is convex\\
	\noalign{\hrule height0.8pt}
	\end{tabular}
\end{center}
\end{table}
\noindent
Our intrinsic topology is the same as their topology. 
However, we are interested in 
the equivalences between the 
intrinsic topologies and the Hilbertian topologies defined 
in section \ref{sec:3}.  

\smallskip
\section{Comparison of Conditions with the CKS-space}
   \label{sec:5}
In this section we will discuss the 
continuity of various operators and wick products.
This matter is beyond the scope of results in 
\cite{ghor}.
   
Let us consider again a CKS-space $[\ce]_{\a}\subset
(L^{2})\subset [\ce]_{\a}^{*}$ associated with a sequence
$\{\a(n)\}$ of positive real numbers. 
Due to the discussion in sections \ref{sec:3} and \ref{sec:4}, 
we conclude that for a CKS-space 
it is reasonable to assume the {\em four essential 
conditions}: (A1), (A2), near-(B2), near-($\wt{\text{B}}2$).

On the other hand, in \cite{kks} the 
following three conditions are imposed in order to prove 
the continuity of various linear operators acting on the 
spaces $[\ce]_{\a}$ and $[\ce]_{\a}^{*}$:

\begin{itemize}
\item[(C1)] There exists a constant $c_{1}$ such that 
for all $n\leq m$,
\begin{equation}
   \a(n) \leq c_{1}^{m} \a(m).  \notag
\end{equation}
\item[(C2)] There exists a constant $c_{2}$ such that 
for all $n$ and $m$,
\begin{equation}
  \a(n+m) \leq c_{2}^{n+m} \a(n) \a(m).  \notag
\end{equation}
\item[(C3)] There exists a constant $c_{3}$ such that 
for all $n$ and $m$,
\begin{equation}
\a(n) \a(m) \leq c_{3}^{n+m} \a(n+m).  \notag
\end{equation}
\end{itemize}

\noindent
It is shown in \cite{kks} that (C3) implies (C1). In the
next two theorems we will show that conditions near-(B2) 
and near-($\wt{\text{B}}2$) imply conditions (C2) and
(C3), respectively.

\begin{theorem} \label{thm:k}
If a sequence $\{\a(n)\}$ of positive real numbers 
satisfies condition near-(B2) and $\a(0)\geq 1$, then it 
satisfies condition (C2).
\end{theorem}

\begin{pf}
Since $\{\a(n)\}$ satisfies condition near-(B2), it is
equivalent to a sequence $\{\l(n)\}$ of positive real 
numbers such that $\{\l(n)/n!\}$ is log-concave. Apply 
Equation (\ref{eq:add}) to the sequence $\b(n) = 
\l(n)/\l(0)$.  Then we get
\begin{equation} \label{eq:k-1}
 \l(n+m) \leq \l(0)^{-1} 2^{n+m} \l(n) \l(m), \qquad
 \forall n, m\geq 0.
\end{equation}
On the other hand, recall that $\{\a(n)\}$ and $\{\l(n)\}$ 
are equivalent. Hence there exist constants $K_{1}, K_{2}, 
c_{1}, c_{2}>0$ such that
\begin{equation} \label{eq:k-2}
 K_{1}c_{1}^{n}\l(n) \leq \a(n) \leq K_{2}c_{2}^{n}\l(n).
\end{equation}
From Equations (\ref{eq:k-1}) and (\ref{eq:k-2}) we can 
easily derive that
\begin{equation}
 \a(n+m) \leq \l(0)^{-1} K_{1}^{-2}K_{2}\big(2c_{1}^{-1}
 c_{2}\big)^{n+m} \a(n)\a(m), \qquad \forall n, m\geq 0.
   \notag
\end{equation}
Let $c=\max\big\{1, \,\l(0)^{-1} K_{1}^{-2}K_{2},\, 
2c_{1}^{-1} c_{2}\big\}$. Then the last inequality 
implies that
\begin{equation}
 \a(n+m) \leq c^{2(n+m)} \a(n) \a(m), \qquad
  \forall n+m\geq 1.  \notag
\end{equation}
However by assumption $\a(0)\geq 1$ and so this inequality also
holds for $n=m=0$. Thus the sequence $\{\a(n)\}$ satisfies
condition (C2).
\end{pf}

\begin{theorem}
If a sequence $\{\a(n)\}$ of positive real numbers 
satisfies condition near-($\wt{\text{B}}2$) and $0<\a(0)
\leq 1$, then it satisfies condition (C3).
\end{theorem}

\begin{pf}
Since $\{\a(n)\}$ satisfies condition 
near-($\wt{\text{B}}2$), it is equivalent to a sequence 
$\{\l(n)\}$ of positive real numbers such that $\{{1\over
n!\l(n)}\}$ is log-concave. Apply Equation (\ref{eq:add}) 
to the sequence $\b(n) = \l(0)/\l(n)$. 
Then we get
$$
	\l(n)\l(m) \leq \l(0) 2^{n+m} \l(n+m)
	\qquad, \forall n, m\geq 0.
$$
Then we
can repeat similar arguments as in the proof of Theorem
\ref{thm:k} to show that the sequence $\{\a(n)\}$ 
satisfies condition (C3).
\end{pf}

For the rest of this section we assume that $u\in 
C_{+, 1/2}$ satisfies conditions (U0) (U2) (U3). We will
state several theorems concerning various continuous linear 
operators acting on $[\ce]_{u}$ and $[\ce]_{u}^{*}$. These 
theorems follow from section 3 of the paper \cite{kks} as 
a consequence of the above Theorem \ref{thm:5-1}. However, 
we point out that they can be proved independently without 
using the corresponding results in the paper \cite{kks}.

The next theorem corresponds to Theorem 3.1 in \cite{kks}.

\begin{theorem} \label{thm:do}
For any $y\in\ce'$, the differential operator $D_{y}$ is
a continuous linear operator from $[\ce]_{u}$ into itself.
\end{theorem}

The next theorem corresponds to Theorem 3.2 in \cite{kks}.

\begin{theorem}
For any $y\in\ce'$, the translation operator $T_{y}$ is
a continuous linear operator from $[\ce]_{u}$ into itself.
\end{theorem}

The next theorem corresponds to a fact on page 323 in 
\cite{kks}.

\begin{theorem}
For any $z\in\spc$, the scaling operator $S_{z}$ is a
continuous linear operator from $[\ce]_{u}$ into itself.
\end{theorem}

For $a, b\in \spc$, define the Fourier-Gauss transform
$\cg_{a, b}\vf$ of $\vf\in [\ce]_{u}$ by
\begin{equation} \label{eq:fg}
 \cg_{a, b}\vf(x) = \int_{\ce'} \vf(ay+bx)\,d\m(y). 
\end{equation}

\begin{theorem} \label{thm:fg}
For any $a, b\in\spc$, the Fourier-Gauss transform 
operator $\cg_{a, b}$ is a continuous linear operator 
from $[\ce]_{u}$ into itself.
\end{theorem}

For those operators in Theorems \ref{thm:do} to 
\ref{thm:fg}, their adjoints are continuous linear 
operators from $[\ce]_{u}^{*}$ into itself. All properties 
regarding to these operators in the book \cite{kuo96} are 
all valid with suitable modification. In particular, the 
integral kernel operators in Chapter 10 and white noise
integration in Chapter 13 can be extended to the Gel'fand 
triple $[\ce]_{u}\subset (L^{2}) \subset [\ce]_{u}^{*}$.

\begin{theorem}
The space $[\ce]_{u}^{*}$ is closed under the Wick 
product and the mapping $(\F, \Q) \mapsto \F\dm\Q$
is jointly continuous from $[\ce]_{u}^{*}\times [\ce]_{u}^{*}$
into $[\ce]_{u}^{*}$ with respect to the inductive limit 
convex topology.
\end{theorem}

\begin{pf}
By Fact \ref{fact:2-6}, $u^{*}$ belongs to $C_{+, 1/2}$ 
and is increasing and (log, $x^{2}$)-convex. 
Hence we may apply the Remark of Lemma \ref{lem:4-1} to $u^{*}$ 
to see that $u^{*}$ and $(u^{*})^{2}$ are equivalent. 
Hence by Theorem \ref{thm:3-1} we can see that 
the space $S[\ce]_{u}^{*}$ is closed under multiplication 
and hence $[\ce]_{u}^{*}$ is closed under Wick product by 
definition. Similarly to Theorem 3.5 in \cite{kks}, we can 
show the inequality 
$$
\|\F\dm\Q\|_{-q, (u)} \leq c \, \|\F\|_{-p, (u)}
 \|\Q\|_{-p, (u)}, \quad \forall 
 \F, \Q\in [\ce_{p}]_{u}^{*}
$$
for any $p$, $c > 1$ and $q > p + \gamma(c)$ with a 
suitable constant $\gamma(c)$.  
The joint continuity can be proved 
applying Lemma A in Appendix 
to $X' = [\ce]_{u}^{*}$. 
\end{pf}

The next theorem is for the Wick product of test functions.
It corresponds to Theorem 3.4 in \cite{kks}.
\begin{theorem}
The space $[\ce]_{u}$ is closed under the Wick 
product and the mapping $(\vf, \q) \mapsto 
\vf\dm\q$ is continuous from $[\ce]_{u}\times [\ce]_{u}$
into $[\ce]_{u}$.
\end{theorem}

For the pointwise multiplication of test functions we have
the next theorem which corresponds to a fact on page 326 
in \cite{kks}.

\begin{theorem} \label{thm:pm}
The space $[\ce]_{u}$ is closed under pointwise 
multiplication and the mapping $(\vf, \q) \mapsto 
\vf\q$ is continuous from $[\ce]_{u}\times [\ce]_{u}$
into $[\ce]_{u}$.
\end{theorem}

\smallskip
\section*{Appendix}
\appendix

To prove the following Lemma A
we borrow the idea in parts from the proof of Lemma 2.1
in \cite{ky89}.
Let us give a short remark. 
It is shown in \cite{ky89} 
that the mapping from a dual space $\ce^{*}$ of $\ce$
to a n-fold symmetirc tensor product space 
$(\ce^{\otimes n})^{*}_{symm}$ is continuous
with respect to 
the inductive limit convex topology.
The nuclearity of $\ce$ 
plays an essential role to prove this fact.
On the other hand, as illustrated in Lemma A, 
the nuclearity of $\ce$ is not an intrinsic assumption 
to verify the continuity 
of the mapping from a direct product space 
$[\ce]_u^{*}\times[\ce]_u^{*}$ to $[\ce]_u^{*}$.
That is, we have 
\smallskip

\noindent
{\bf Lemma A.} 
{\it Let $X$ be a complete $\sigma$-normed space 
with norms $\|\cdot\|_1 \leq \|\cdot\|_2 \leq \cdots \leq 
\|\cdot\|_p \leq \cdots$, and let $X'$ be its dual and 
$\|\cdot\|_{-p}$ the dual norm of $\|\cdot\|_p$.  
Suppose that $X'$ is an algebra with multiplication $xy \in X'$ of 
$x, y \in X'$ and that for any $p\geq 1$ there 
exist an integer $\gamma(p)\geq p$ and a positive constant 
$C(p)$ such that
$$
\|xy\|_{-\gamma(p)} \leq C(p) \|x\|_{-p}\|y\|_{-p}
$$
holds for any $x, y \in X'$ with $\|x\|_{-p},\|y\|_{-p}
< \infty$.
Then  the mapping $(x,y)\mapsto xy$
is jointly continuous from $X'\times X'$ 
into $X'$ with respect to the inductive limit convex 
topology}. 
\begin{pf}
Recall that fundamental neighborhoods of $x$ 
can be given in the form
$$
V(x;q,\{\epsilon_p\}_{p\geq q}) 
=\hbox{conv}\left(\bigcup_{p\geq q}\{z; \|z\|_{-p} 
< \epsilon_p\}\right) + x,
$$
for a given $q \geq 1$ and a positive sequence 
$\{\epsilon_p\}_{p \geq q}$ (cf. \cite{ky89}). 
Here "conv" means 
the convex hull, just the collection of all finite convex sums. 
For $V(0;q,\{\epsilon_p\}_{p\geq q})$, 
put 
$$
\delta_p = \min\left\{1, {\epsilon_{\gamma(p)} \over C(p)}\right\}
\quad \hbox{for } \ p \geq q.
$$
For $x, y \in V(0;q,\{\delta_p\}_{p\geq q})$, there 
exist $\{\alpha_p\}_{p\geq q}^N, \{\beta_p\}_{p\geq q}^{N'}$ 
and $\{x_p\}_{p\geq q}^N, \{y_p\}_{p\geq q}^{N'}$ 
such that 
$$
x = \sum_{p=q}^N \alpha_p x_p, \ \
y = \sum_{p=q}^{N'} \beta_p y_p \quad \hbox{and} \quad
\|x_p\|_{-p} < \delta_p, \ \ 
\|y_p\|_{-p} < \delta_p, 
$$
$$
\sum_{p=q}^{N'} \alpha_p = 1, \ \
\sum_{p=q}^{N'} \beta_p = 1, \ \ 
\alpha_p \geq 0, \ \ \beta_p \geq 0.
$$
Then put $\ell_0 = \min \{\gamma(p) : p \geq q\}$ and put 
$$
z_\ell = {1 \over \lambda(\ell)} \sum_{\gamma(k)=\ell} \sum_{p\vee p' = k} 
\alpha_p\beta_{p'} x_py_{p'}, \quad
\lambda(\ell) = \sum_{\gamma(k)=\ell} \sum_{p\vee p' = k} 
\alpha_p\beta_{p'}.
$$
for $\ell \geq \ell_0$.  Now estimate norms of $z_\ell's$. 
Since $\|x_p\|_{-k} \leq \|x_{p}\|_{-p} < \delta_p$ and 
$\|y_{p'}\|_{-k} \leq \|y_{p'}\|_{-p'} < \delta_{p'}$ 
for $k \geq p\vee p'$, we have 
\begin{align*}
\|z_\ell\|_{-\ell}  \leq &
{1 \over \lambda(\ell)} \sum_{\gamma(k)=\ell} \sum_{p\vee p' = k} 
\alpha_p\beta_{p'} \|x_py_{p'}\|_{-\ell} \cr
 \leq & {1 \over \lambda(\ell)} \sum_{\gamma(k)=\ell} \sum_{p\vee p' = k} 
\alpha_p\beta_{p'} C(k)\|x_p\|_{-k}\|y_{p'}\|_{-k} \cr
 \leq & {1 \over \lambda(\ell)} \sum_{\gamma(k)=\ell} \sum_{p\vee p' = k} 
\alpha_p\beta_{p'} C(k)\|x_p\|_{-p}\|y_{p'}\|_{-p'} \cr
 \leq & {1 \over \lambda(\ell)} \sum_{\gamma(k)=\ell} \sum_{p\vee p' = k} 
\alpha_p\beta_{p'} C(k)\delta_{k} \cr
 \leq & \epsilon_\ell. \cr
\end{align*}
Since 
$$
xy = \sum_{\ell\geq \ell_0} \lambda(\ell) z_\ell, \quad 
\sum_{\ell\geq \ell_0} \lambda(\ell) = 1, \quad
\ell_0 \geq q,
$$
we obtain $xy \in V(0;q,\{\epsilon\}_{p\geq q})$.
Hence the product is jointly continuous at $0$.  

Next we show the joint continuity at $(x_0,y_0)$. 
Suppose that $\|x_0\|_{-p_0}, \|y_0\|_{-p'_0} < \infty$.
For any given $V(0;q,\{\epsilon_p\}_{p\geq q})$, put 
$q_0 = \max\{q,p_0, p'_0\}$ and 
$$
\delta_p =  \min\left\{1, {\epsilon_{\gamma(p)} 
\over 3C(p)(1 + \|x_0\|_{-p} + |y_0\|_{-p})}\right\} \quad 
\hbox{for } \ p \geq q_0
$$
and  take their neighborhoods as 
$V(x_0;q_0,\{\delta_p\}_{p\geq q_0})$ 
and $V(y_0;q_0,\{\delta_p\}_{p\geq q_0})$ 
and let $x$ and $y$ be in these neighborhoods, respectively.  
Then we have   
$$
x = x_0 + \sum_{p=q_0}^N \alpha_p x_p, \quad
y = y_0 + \sum_{p=q_0}^{N'} \beta_p y_p 
$$
as above.  Then we see 
$$
xy - x_0y_0 =  (x - x_0)(y - y_0) 
+ (x - x_0)y_0 + (y - y_0)x_0.
$$
The first term of the right hand side belongs to 
$V(0;q,\{{1 \over 3}\epsilon_p\}_{p\geq q})$. Observe the 
second term. Put $\ell_0 = \min \{\gamma(k) ; k \geq q_0\}$ and 
$$
z'_\ell = {1 \over \lambda'(\ell)} \sum_{\gamma(k) = \ell, k \geq q_0} 
\alpha_k x_k y_0, \quad 
\lambda'(\ell) = \sum_{\gamma(k) = \ell, k \geq q_0} \alpha_k.
$$
Then 
\begin{align*}
\|z'_\ell\|_{-\ell}  \leq & 
  {1 \over \lambda'(\ell)} 
  \sum_{\gamma(k) = \ell, k \geq q_0} \alpha_k \|x_k y_0\|_{-\ell }\cr
 \leq & {1 \over \lambda'(\ell)} \sum_{\gamma(k) = \ell, k \geq q_0} 
 \alpha_k C(k)\|x_k\|_{-k}\|y_0\|_{-k} \cr
 \leq & {1 \over \lambda'(\ell)} \sum_{\gamma(k) = \ell, k \geq q_0} 
 \alpha_k C(k)\delta_k \|y_0\|_{-k} < {1 \over 3}\epsilon_\ell
\end{align*}
This implies that $(x-x_0)y_0 \in 
V(0;q,\{{1 \over 3}\epsilon_p\}_{p\geq q})$. 
In the same way, we see $(y-y_0)x_0 \in 
V(0;q,\{{1 \over 3}\epsilon_p\}_{p\geq q})$. 
Thus we have $xy - x_0y_0 \in V(0;q,\{\epsilon_p\}_{p\geq q})$.
We complete the proof. 
\end{pf}

\bigskip

\noindent
{\bf Acknowledgements.} The authors express their hearty 
thanks to the referee for his comments improving this paper. 
N.~Asai wants to thank the Daiko
Foundation and the Kamiyama Foundation for research support. 
I.~Kubo wishes to express thanks to the Ministry of 
Education, Science, Sports and Culture of Japan 
for Grant-in-Aid for Scientific Research 2000 (Nr. 10304006).
H.-H.~Kuo is grateful for financial supports from the 
Academic Frontier in Science (AFS) of Meijo University and 
the Luso-American Foundation. He wants to thank AFS and
CCM, Universidade da Madeira for the warm hospitality during 
his visits (February 15--21, 1998 and March 1--7, 1999 to
AFS, July 22-August 20, 1999 to CCM.) In particular, he gives 
his deepest appreciation to Professors T. Hida and K. Sait\^o
(AFS) and M. de Faria and L. Streit (CCM) for arranging 
the visits. H.-H.~Kuo also wants to thank Professor Y.-J.~Lee
for arranging his visit to Cheng Kung University in the
spring of 1998. At that time this joint research project 
started. He thanks the financial support from the National 
Science Council of Taiwan.


\end{document}